\documentclass[12pt,reqno]{amsart}
\usepackage{amssymb,delarray}
\usepackage{amsfonts}
\usepackage{epsfig}

\makeindex{}

\def\le{\leqslant}
\def\ge{\geqslant}

\newtheorem{thm}{Theorem}[section]
\newtheorem{lem}[thm]{Lemma}
\newtheorem{prop}[thm]{Proposition}

\newtheorem{cor}[thm]{Corollary}
\newtheorem{rem}[thm]{Remark}
\newtheorem{ex}[thm]{Example}

{\catcode`\@=11
\gdef\n@te#1#2{\leavevmode\vadjust{%
 {\setbox\z@\hbox to\z@{\strut#1}%
  \setbox\z@\hbox{\raise\dp\strutbox\box\z@}\ht\z@=\z@\dp\z@=\z@%
  #2\box\z@}}}
\gdef\leftnote#1{\n@te{\hss#1\quad}{}}
\gdef\rightnote#1{\n@te{\quad\kern-\leftskip#1\hss}{\moveright\hsize}}
\gdef\?{\FN@\qumark}
\gdef\qumark{\ifx\next"\DN@"##1"{\leftnote{\rm##1}}\else
 \DN@{\leftnote{\rm??}}\fi{\rm??}\next@}}

\begin{document}
\baselineskip=14.pt plus 2pt 

\title[Real "Dif=Def" problem]{Surfaces with DIF$\ne$DEF real structures}
\author[V.M.~Kharlamov and Vik.S.~Kulikov]{V.M.~Kharlamov and Vik. S.~Kulikov}
\address{
UFR de Math\'ematiques et IRMA \\ Universit\'e Louis Pasteur et
CNRS \\ 7 rue R\'en\'e-Descrates \\ 67084, Strasbourg, Cedex,
France}
\email{ kharlam@math.u-strasbg.fr}
\address{Steklov Mathematical Institute\\
Gubkina str., 8\\
119991 Moscow \\
Russia} \email{kulikov@mi.ras.ru}

\dedicatory{} \subjclass{}
\thanks{This work is done during the stay of the second
author in Strasbourg university.
The first author is a member of Research Training Network RAAG
CT-2001-00271. The second author is partially supported by RFBR
({\rm No.} 05-01-00455) and  NWO-RFBR-047.011.2004.026.
 }
\keywords{}
\begin{abstract}
We study real Campedelli surfaces up to real deformations and
exhibit a number of such surfaces which are equivariantly
diffeomorphic but not real deformation equivalent.
\end{abstract}

\maketitle \setcounter{tocdepth}{2}

\def\st{{\sf st}}

\setcounter{section}{0}
\section*{Introduction}
The real DIF=DEF problem is at least as old as the complex one. As
in the complex DIF=DEF problem it is a question of interaction
between two basic equivalence relations:
by diffeomorphisms of real structures, and by deformations of varieties
together with real structures.

A \emph{real structure} on a complex surface~$X$ is an
anti-holomorphic involution $X\to X$. A complex surface supplied
with a real structure is called a \emph{real surface}. A
\emph{deformation} of
surfaces is a proper
holomorphic submersion $p : Z\to D$, where $Z$ is a
$3$-dimensional complex
variety and $D\subset\mathbb C$ is a disk. If $Z$ is real and $p$
is equivariant, the deformation is called real. Two real
surfaces~$X'$ and~$X''$ are called \emph{deformation equivalent}
if they can be connected by a chain $X'=X_0$, \dots, $X_k=X''$ so
that $X_i$ and $X_{i-1}$ are isomorphic to real fibers of a real
deformation.

Under these definitions, {\it up to a diffeomorphism the real
structure is preserved under deformation.} So the problem is in
what extent the diffeomorphic type of the real structure
determines the deformation type. Namely, let call a real surface
$X$ to be {\it quasi-simple} if it is deformation equivalent to
any other real surface $X'$ such that, first, $X'$ is deformation
equivalent to $X$ as a complex surface, and, second, the real
structure of $X'$ is diffeomorphic to the real structure of $X$.
Thus, we understand the real DIF=DEF problem as the question {\it
are there non quasi-simple real surfaces}? (Note that in the case
of curves the response to such a question is in negative: {\it any real
curve is quasi-simple}.)

The first quasi-simplicity result
belongs to F.~Klein and\,
L.~Schl\"{a}fli \cite{Kl} and concerns real cubic surfaces
in the projective $3$-space.
In fact, the quasi-simplicity holds for many other special classes
of surfaces. It is observed for rational surfaces (A.~Degtyarev and
V.~Kharlamov~\cite{DKrat}), for real Abelian surfaces
(follows from A.~Comessatti~\cite{Com-abel}), for geometrically
ruled real surfaces (J.-Y.~Welschinger~\cite{W}), for real
hyperelliptic surfaces (F.~Catan\-ese and P.~Frediani~\cite{CF}),
for real $K3$-surfaces (follows from V.~Nikulin~\cite{NIK}), and
for real Enriques surfaces (A.~Degtyarev and V.~Kharlamov; the
quasi-simplicity statement was announced in~\cite{DK}, and the
complete list of deformation classes of real Enriques surfaces was
obtained in collaboration with I.~Itenberg in~\cite{Book}; note
also that quasi-simplicity of hyperelliptic and Enriques surfaces
extends to quasi-simplicity of the quotients of Abelian and
$K3$-surfaces by certain finite group actions, see \cite{Duke}).

Whether elliptic surfaces and irrational ruled surfaces
quasi-simple is, as far as we know, still an open question.

It was natural to expect that such a simple behaviour would no
longer take place for more complicated surfaces, like those of
general type. However, probably because of lack of convenient
deformation invariants not covered by the differential topology of
the real structure, no any example of non quasi-simple real
surfaces (or real varieties of higher dimension) was known. The
main result of this paper is providing such examples. Namely, we
prove that {\emph{ the Campedelli surfaces {\rm (see the
definition in Section \ref{campcovering})} are not quasi-simple: there
exist real Campedelli sufaces which have diffeomorphic real
structures without being deformation equivalent}.

Let us notice that existence of non quasi-simple families of
surfaces of general type does not prevent certain particular
classes of surfaces of general type from being quasi-simple. And
examples of quasi-simple real surfaces of general type do exist.
One such example is given by real Bogomolov-Miyaoka-Yau surfaces,
that is, surfaces covered by a ball in $\mathbb C^2$, see
\cite{KK}. In fact, in \cite{KK} it is also shown that there exist
diffeomorphic, in fact complex conjugated, Bogomolov-Miyaoka-Yau
surfaces which are not real and thus, being rigid, they are not
deformation equivalent. These surfaces are counter-examples to the
$\rm Diff=Deff$ problem in complex geometry.

The first counter-examples to the $\rm Diff=Deff$ problem in the
complex geometry of surfaces belong to Manetti~\cite{Man}. They
are not involving the complex conjugation. Already their existence
explains why we need to fix complex deformation class in the
definition of quasi-simplicity of real varieties. Moreover, our
examples of diffeomorphic but not deformation equivalent real
structures are closely related to Manetti's examples. In fact, to
establish a diffeomorphism we follow Manetti's approach, and to
study the deformation equivalence we use the full description of
the Campedelli surfaces given by Miayoka \cite{Mi}.

The paper is organized as follows. In Section 1, we collect
essentially known results on complex Campedelli surfaces adapting
them to our needs and making emphasis on representing Campedelli
surfaces as Galois coverings of $\mathbb P^2$. In Section 2, we
begin our study of real structures on Campedelli surfaces and give
a kind of classification of real structures on such surfaces.
Section 3 is devoted to a study of real structures up to
diffeomorphisms and up to deformations. In Section 4, we apply the
technique developed to construct real surfaces which have
diffeomorphic real structures without begin deformation
equivalent. Related remarks are collected in Section 5.

\section{Moduli space of Campedelli surfaces}

\subsection{Campedelli surfaces as branched Galois
coverings of the projective plane}\label{campcovering} Let $X$ be
{\it a Campedelli surface}, that is a minimal surface of general
type which has $p_g=q=0, K_X^2=2$, and $\pi_1(X)=(\mathbb
Z/2\mathbb Z)^3$. Denote by $X_{can}={\rm Proj}(\sum_m H^0(X;
mK))$ the canonical model of $X$, by $\tilde X$ the universal
covering of $X$, by $G_{un}$ the Galois group of this universal
covering, and by $\tilde X_{can}$ the canonical model of $\tilde
X$. Note that $\tilde X_{can}$ and $X_{can}$ have at most simple
double points as singularities, so that $\tilde X_{can}$ is the
universal covering of $X_{can}$. The universal coverings $\tilde
X\to X$ and $\tilde X_{can}\to X_{can}$ have the same Galois
group, so that $X_{can}=\tilde X_{can}/G_{un}$.

According to \cite{Mi}, Theorem 9, the following statement holds.

\begin{thm} \label{thmMi}
The canonical map imbeds $\tilde X_{can}$ in $\mathbb P^6$. With
respect to suitable homogeneous coordinates $w_0, \dots, w_6$ in
$\mathbb P^6$, this image of $\tilde X_{can}$ is given by
equations
\begin{equation}
w_i^2=a_iw_0^2 +b_iw_1^2+c_iw_2^2,\quad  a_i, b_i, c_i \in\mathbb
C, \qquad i=3,4,5,6, \label{quad}
\end{equation}
and the group $G_{un}=(\mathbb Z/2\mathbb Z)^3$ acts on $\tilde
X_{can}$ by diagonal projective transformations:
 $g^*(w_j)=\pm w_j$ for any $g\in G_{un}$.\qed
\end{thm}

As equations (\ref{quad}) and Theorem \ref{thmMi} show, the whole
group $\widetilde G \simeq (\mathbb Z/2\mathbb Z)^6\subset PGL
(6,\mathbb C)$ of diagonal involutions ($g^*(w_j)=\pm w_j$ for any
$g\in \tilde G$) acts on $\tilde X_{can}$ and the following
statement holds.

\begin{cor} \label{claim}
The quotient space $\tilde X_{can}/\widetilde G$ is isomorphic to
$\mathbb P^2$ and the quotient map $\tilde X_{can}\to \tilde
X_{can}/\widetilde G$ is a Galois covering of $\mathbb P^2$ with
Galois group $\widetilde G \simeq (\mathbb Z/2\mathbb Z)^6$
branched along seven lines given by equations
$$ \begin{array}{l}
z_i=0, \qquad i=0,1,2, \\
z_i=  a_iz_0 +b_iz_1+c_iz_2,\qquad i=3,4,5,6,
\end{array}$$ where
$z_0,z_1,z_2$ are homogeneous coordinates in $P^2=\tilde
X_{can}/\widetilde G$.

The canonical model $X_{can}$ of $X$ is a Galois covering of
$\mathbb P^2$ with Galois group $G\simeq (\mathbb Z/2\mathbb Z)^3$
branched along the same lines.\qed
\end{cor}

Let us underline that in the above statements the choice of the
equations and the coverings is not arbitrary.

\subsection{Few basic facts on Galois coverings}\label{basic}
Recall that a Galois covering of a smooth algebraic variety $Y$ is
a finite morphism $h:X\to Y$ of a normal algebraic variety $X$ to
$Y$ such that the function fields imbedding $\mathbb C (Y)\subset
\mathbb C (X)$ induced by $h$ is a Galois extension. As is well
known, a finite morphism $h:X\to Y$ is a Galois covering with
Galois group $G$ if and only if $G$ coincides with the group of
covering transformations and the latter acts transitively on every
fiber of $h$. Besides, a finite branched covering is Galois if and
only if the unramified  part of the covering (i.e., the
restriction to the complements of the ramification and branch
loci) is Galois. In addition, a branched covering is determined up
to isomorphism by its unramified part. Moreover, a covering map
from the unramified part $U_1\subset X_1$ of one branched
covering, $h_1:X_1\to Y_1$, to the unramified part $U_2\subset
X_2$ of another one, $h_2:X_2\to Y_2$, induces a covering morphism
$X_1\to X_2$ between these branched coverings if the extension of
the morphism of underlying unbranched varieties, $h_1(U_1)\subset
Y_1$ and $h_2(U_2)\subset Y_2$, to the branch loci is given. Let
us recall also that an unramified covering is Galois with Galois
group $G$ if and only if it is a covering associated with an
epimorphism of the fundamental group of the underlying variety to
$G$, and, in particular, the Galois coverings with abelian Galois
group $G$ are in one-to-one correspondence with epimorphisms to
$G$ of the first homology group with integral coefficients. All
these results are well known and their most nontrivial part can be
deduced, for example, from the Grauert-Remmert existence theorem
\cite{G-R}.

In what follows we deal with Galois coverings with Galois group
$G\simeq (\mathbb Z /2\mathbb Z )^k$. Galois groups are considered
up to isomorphism, and two Galois coverings $h_1 : X_1\to Y$ and
$h_2 : X_2\to Y$ with Galois groups $G_1$ and $G_2$ are said to be
{\it equivalent} if there exist a biregular map $f: X_1\to X_2$
and an isomorphism $F : G_1\to G_2$ such that $h_2\circ f=h_1$ and
$F(g)f(x)=f(gx)$ for any $x \in X_1$ and $g \in G_1$.

\subsection{Galois coverings of $\mathbb P^2$ with Galois group
$(\mathbb Z/2\mathbb Z)^k$ branched along seven lines}\label{linecovering}
Let
$\mathcal L=L_0\cup \dots \cup L_6$ be an arrangement of seven
distinct numbered lines in $P^2$. The simple loops $\lambda _i,
0\le i\le 6,$ around the lines $L_i$ generate $H_1(\mathbb P ^2
\setminus \mathcal L ,\mathbb Z )\simeq \mathbb Z ^{6}$. They are
subject to the relation
$$\lambda _0+\dots +\lambda _6=0.$$

The natural epimorphism $\widetilde \varphi : H_1(\mathbb P ^2
\setminus \mathcal L, \mathbb Z )\to H_1(\mathbb P ^2 \setminus
\mathcal L,\mathbb Z/2\mathbb Z )\simeq (\mathbb Z /2\mathbb Z )
^6$ defines a particular Galois covering of $P^2$ ramified in
$\mathcal L$. We call it universal and denote by $\widetilde g:
\widetilde Y\to \mathbb P ^2$. The following statement, which is a
straightforward consequence of the general results on branched
coverings mentioned in Section \ref{basic}, precises, in
particular, at what sense it is universal.

\begin{prop}\label{0.1}
Galois coverings with Galois group $G\simeq (\mathbb Z /2\mathbb Z
)^k$ branched along $\mathcal L$ exist if and only if $k\le 6$.
Their equivalence classes are in one-to-one correspondence with
epimorphisms $H_1(\mathbb P^2\setminus\mathcal L)\to G$ considered
up to automorphisms of $G$. If $g: Y \to \mathbb P ^2$ is a Galois
covering with Galois group $G\simeq (\mathbb Z /2\mathbb Z )^k$
branched along $\mathcal L$, then there exists a Galois covering
$h: \widetilde Y\to Y$ such that $\widetilde g=g\circ h$. \qed
\end{prop}

Without loss of generality, we can assume that the universal
Galois covering $\widetilde g: \widetilde Y\to \mathbb P ^2$ is
associated with the epimorphism $\widetilde \varphi : H_1(\mathbb
P ^2 \setminus \mathcal L , \mathbb Z )\to (\mathbb Z /2\mathbb Z
)^{6}$ sending $\lambda _0$ to $(1,\dots , 1)$ and $\lambda _i$
with $1\leq i\leq 6$ to $(0,\dots,0,1,0,\dots,0)$ with $1$ in the
$i$-th place.

Let $(v_1,v_2)$ be affine coordinates in $\mathbb C ^2=\mathbb P
^2\setminus L_0$ and $l_i(v_1,v_2)=0$, $1\le i\le 6,$ be a linear
equation of $L_i \cap \mathbb C ^2$. The function field
$K_{u}=\mathbb C (\widetilde Y)$ of $\widetilde Y$ is the abelian
extension $\mathbb C (\widetilde Y)=\mathbb C (v_1,v_2,w_1,\dots
,w_{6})$ of the function field $K=\mathbb C (v_1,v_2)$ of $\mathbb
P ^2$ of degree $2^{6}$ determined by $w_i^2=l_i$, $i=1,\dots, 6$.
(In other words, the pull-back of $\mathbb P^2\setminus L_0$ in
$\widetilde Y$ is naturally isomorphic to the normalization of the
affine subvariety of $\mathbb C^{8}$ given in affine coordinates
$v_1,v_2,w_1,\dots, w_{6}$ by equations $w_1^2=l_1,\dots ,
w_{6}^2=l_{6}$.)

The action of $\gamma= (\gamma_1,\dots ,\gamma _{6})\in (\mathbb Z
/2\mathbb Z )^{6}$ on $K_{u}$ is given by
$$
\gamma (w^a) =(-1)^{( \gamma, a)} w^a,
$$
where for any multi-index $a=(a _1,\dots ,a_{6})$, $0\leq a _i\leq
1$, we put
$$
w^a =\prod_{i=1}^{6}w_i^{a _i}.$$ Therefore,
$\mbox{Gal}(K_u/\mathbb C(v_1,v_2))= (\mathbb Z /2\mathbb Z )^{6}$
and
$$
K_{u}=\bigoplus_{ 0\le a _i\le 1} \mathbb C (v_1,v_2)w^a
$$
is a decomposition of the vector space $K_{u}$ over $\mathbb C
(v_1,v_2)$ into a finite direct sum of degree $1$ representations
of $(\mathbb Z /2\mathbb Z )^{6}$.

Let $\varphi : H_1(\mathbb P^2 \setminus \mathcal L,\mathbb Z )\to
(\mathbb Z /2\mathbb Z)^k$ be an epimorphism given by $\varphi
(\lambda _i)=(a _{i,1},\dots ,a _{i,k})$, where $a _{0,j}+\dots +a
_{6,j}\equiv 0\, \text{mod}\, 2$ for every $j=1,\dots , k$, and
let $g:Y \to \mathbb P ^2$ be the Galois covering associated with
$\phi$. This covering is ramified in the union $\mathcal
L^\varphi$ of lines $L_i\subset\mathcal L$ with
$\varphi(\lambda_i)\ne 0$. The epimorphism $\varphi$ factors
through a unique epimorphism $\psi : (\mathbb Z /2\mathbb Z
)^{6}\to (\mathbb Z /2\mathbb Z)^{k}$, so that, by Proposition
\ref{0.1}, the covering $g$ factors through a unique Galois
covering $h: \widetilde Y\to Y$. The latter determines the
inclusion $h^{*} :\mathbb C (Y)\to K_{u}$ of the function field
$\mathbb C (Y)$ of $Y$ into the function field $K_{u}=\mathbb
C(\widetilde Y)$. Clearly, $\mbox{Gal}(K_{u}/h^*(\mathbb C
(Y)))=\ker \psi$, the field $h^*(\mathbb C(Y))$ coincides with the
subfield $K_{\varphi}=\mathbb C (v_1,v_2,u_1,\dots ,u_{k})$ of
$K_{u}$, where
\begin{equation} \label{new} u_j=w_1^{a _{1,j}}\cdot .\, .\, . \cdot w_{6}^{a_{6,j}},
\end{equation} and
$$\mbox{Gal}(K_{u}/K_{\varphi})= \{ \,
(\gamma_1,\dots ,\gamma _{6})\in
(\mathbb Z /2\mathbb Z )^{6}\, \mid \, \sum_{i=1}^{6}
a_{i,j}\gamma_{i}\equiv 0\,
(2),\, 1\leq j \leq k\, \}.$$

\subsection{ Resolution of singularities of $Y$}\label{Yresolution} By construction,
$Y$ is a normal surface with isolated singularities. The singular
points of $Y$ can appear only over an $r$-fold point of $\mathcal
L^\varphi$ with $r\ge 2$, i.e., over a point belonging to exactly
$r$ lines $L_{i_1},\dots ,L_{i_r}\in\mathcal L$ with
$\varphi(\lambda_{i_k})\ne 0$, $1\le k\le r$.

\begin{lem} \label{lem 1.1}
{\rm (see, f.e., (\cite{Kh-Ku})} If $p=L_{i_1}\cap L_{i_2}$ is a
$2$-fold point of $\mathcal L^\varphi$ and $\varphi
(\lambda_{i_1}) \ne \varphi (\lambda_{i_2})$, then $Y$ is
non-singular at each point of $g^{-1}(p)$. \qed
\end{lem}

We say that an $r$-fold point $p_{i_1,\dots,i_r}$ of $\mathcal
L^\phi$ is {\it a non-branch point with respect to $\varphi$} if
$\sum_{j=1}^r\varphi(\lambda _{i_j})=0$.

To resolve the singularities of $Y$ we start from a suitable
blow-up of $\mathbb P^2$. First, we blow up all the $2$-fold
non-branch points and all the $r$-fold points of $\mathcal
L^\varphi$ with $r\geq 3$. Second, for each pair
$(p_{i_1,\dots,i_r}, k)$ such that $p_{i_1,\dots,i_r}$ is an
$r$-fold point of $\mathcal L^\varphi$ and
$\sum_{j=1}^r\varphi(\lambda _{i_j})=\varphi(\lambda_{i_k})$, we
effectuate a blow-up with center at the intersection point of the
strict transform of $L_{i_k}$ with the exceptional divisor
$E_{i_1,\dots,i_r}$ blown-up over $p_{i_1,\dots,i_r}$ at the first
series of the blow-ups. The resulting combination of the blow-ups
is denoted by $\sigma :\hat {\mathbb P }^2 \to \mathbb P ^2$.

By $L'_i\subset \hat{\mathbb P ^2}$ we denote the strict transform
of $L_i$, by $E'_p$ with $p=p_{i_1,\dots,i_r}$ the strict
transform of $E_{i_1,\dots,i_r}$, by $E_{p,\,i_k}$ the exceptional
curves of the second series of the blow-ups, and by
$\varepsilon_{p},\varepsilon_{p,\,i_k}\in H_1(\hat{\mathbb P
^2}\setminus \sigma ^{-1}(\mathcal L^\varphi),\mathbb Z )
=H_1(\mathbb P ^2 \setminus \mathcal L^\varphi,\mathbb Z )$ simple
loops around $E'_{p}$ and $E_{p,\,i_k}$, respectively.

The identification $H_1(\hat{\mathbb P ^2} \setminus \sigma
^{-1}(\mathcal L^\varphi),\mathbb Z ) =H_1(\mathbb P ^2 \setminus
\mathcal L^\varphi ,\mathbb Z )$ composed with $\varphi$ provides
an epimorphism $\hat\varphi : H_1(\hat{\mathbb P ^2} \setminus
\sigma ^{-1}(\mathcal L^\varphi),\mathbb Z )\to (\mathbb Z
/2\mathbb Z )^k$. Let consider the associated Galois covering $f:X
\to \hat{\mathbb P ^2}$.

\begin{lem} \label{lem 1.2}  {\rm (\cite{Kh-Ku})}
Let $p=L_{i_1}\cap \dots \cap L_{i_r}$ be an $r$-fold point of
$\mathcal L^\varphi$, $r\ge 2$. Then,
\begin{itemize}
\item[($i$)] $\varepsilon_{p}=\lambda _{i_1}+\dots + \lambda
_{i_r}$, \item[($ii$)] $\varepsilon_{p,\,i_k}=\lambda
_{i_k}+\displaystyle \sum_{j=1}^r\lambda _{i_j}$, \item[($iii$)]
$\varphi(\varepsilon_{p,\,i_k})=0$. \end{itemize} \qed
\end{lem}

 The following theorem is a straightforward consequence of Lemmas
\ref{lem 1.1} and \ref{lem 1.2}.
\begin{thm}\label{resolution}
The Galois coverings $f$ and $g$ are included in the commutative
diagram

\begin{picture}(280,90)
\put(130,70){$X$} \put(160,77){$\nu $}
\put(143,75){\vector(1,0){40}} \put(192,70){$Y$}
\put(133,65){\vector(0,-1){30}}\put(125,47){$f$}
\put(195,65){\vector(0,-1){30}}\put(198,49){$g$}
\put(130,20){$\hat{\mathbb P ^2}$} \put(160,17){$\sigma $}
\put(145,25){\vector(1,0){40}} \put(191,20){$\mathbb P ^2$}
\end{picture}
\newline in which $\nu:X\to Y$ is a resolution of singularities of $Y$.
\qed \end{thm}

\begin{lem}\label{A-D-E}
Suppose that the Galois group of the covering $Y\to\mathbb P^2$ is
$(\mathbb Z /2\mathbb Z)^3$.
Then, a point $q\in Y$ situated over an $r$-fold point $p=p_{i_1,\dots,i_r}=
L_{i_1}\cap\dots\cap L_{i_r}$ of
$\mathcal L^\phi$ is not a canonical singular
point (that is, $q$ is not an $A$-$D$-$E$-singularity) if,
and only if,
either $r> 3$
or $r=3$ and $p$ is not a branch point of $\varphi$.
\end{lem}

\proof To determine the type of a singular point we look
at its resolution provided by $\nu : X\to Y$, see Theorem \ref{resolution}.

If $r=2$ and $\varphi(\lambda_1)\neq \varphi(\lambda_2)$,
then, by Lemma \ref{lem 1.1}, each point $q\in g^{-1}(p)$
is a nonsingular point of $Y$. If, by contrary,
$\varphi(\lambda_1)= \varphi(\lambda_2)$, then the covering
$f:X\to\hat{\mathbb P^2}$ is not branched at $E'_p$
and it splits over $E'_p$ into four copies
of a Galois double covering of $\mathbb P^1$ branched at two
points, so that each of the four points $q\in g^{-1}(p)$
is replaced in the resolution by a rational curve with
self-intersection number $\frac{(-1)\cdot\, 8}{4}=-2$. Hence, in
this case all the four points are of type $A_1$.

If $r=3$ and $p$ is a non-branch point,
then up to a coordinate change in $G$ we have
$\varphi(\lambda_{i_1})=(1,0,0)$,
$\varphi(\lambda_{i_2})=(0,1,0)$, and
$\varphi(\lambda_{i_3})=(1,1,0)$. Therefore,
$f^{-1}(E'_p)$  is a disjoint union of two
rational curves $C_1$ and $C_2$ with self-intersection
$\frac{(-1)\cdot\, 8}{2}=-4$. Hence, the singular points $q\in g^{-1}(p)$ are
not canonical.

Now, let suppose that $r=3$, $p$ is a branch point, and
$\varphi(\lambda_{i_1})$, $\varphi(\lambda_{i_2})$,
$\varphi(\lambda_{i_3})$ are pairwise distinct (note that for a
branch point the latter assumption is equivalent to
$\sum_{j=1}^3\varphi(\lambda _{i_j})\ne\varphi(\lambda_{i_k})$ for
any $1\le k\le 3$). Then, after a coordinate change in $G$ we may
suppose that $\varphi(\lambda_{i_1})=(1,0,0)$,
$\varphi(\lambda_{i_2})=(0,1,0)$,
$\varphi(\lambda_{i_3})=(0,0,1)$. Therefore, over $E'_p$ we get a
Galois covering over $\mathbb P^1$ with Galois group $(\mathbb Z
/2\mathbb Z)^2$ and three branched points, so that $f^{-1}(E'_p)$
is a rational curve with self-intersection $\frac{(-1)\cdot\,
8}{4}=-2$, and, hence, the singular point $q=g^{-1}(p)$ is of type
$A_1$.

Next, let treat the case when $r\ge 3$ and there is at least one
$k$ such that $\sum_{j=1}^r\varphi(\lambda
_{i_j})=\varphi(\lambda_{i_k})$. Then: $p$ is a branch point,
$\sigma^{-1}(p)=E'_p +\sum^s_{j=1} E'_{p,k_j}$ where
$(E'_p)^2=-(s+1)$ and $(E'_{p,k_1})^2=\dots=(E'_{p,k_s})^2=-1$;
$E'_p$ is a branch curve of $f$, but $E'_{p,k_1},\dots,
E'_{p,k_s}$ are not branch curves of $f$. Therefore, each of
$f^{-1}(E_{p,k_j}), 1\le j\le s$, splits into a disjoint union of
four $(-2)$-curves, while $f^{*}(E'_p)=2C_1+\dots +2C_{2^n}$,
where $2^n$ is the index in $G$ of the subgroup
$G_{i_1,\dots,\,i_r}$ generated by $\varphi(\lambda_{i_1}), \dots
,\varphi(\lambda_{i_r})$ and $C_1, \dots, C_{2^n}$ are copies of a
Galois covering of $E'_p$ of degree $2^{2-n}$ (recall that $\deg
f=8$) branched at $r-s$ points. Thus, for each $i=1,\dots,2^n$ we
have
$$\begin{array}{l}
(C_i^2)_{X}=-2^{1-n}(s+1), \\
g(C_i)=2^{-n}(r-s)-2^{2-n} +1=2^{-n}(r-s-4)+1,\end{array}
$$
where $0\leq n\leq 2$. If $g^{-1}(p)$ consists of canonical
singularities, then $(C_i^2)_{X}=-2$ and $g(C_i)=0$, that is,
$$\begin{array}{l} 2^{1-n}(s+1)=2, \\
2^{-n}(r-s-4)+1=0.\end{array} $$ The only solutions are $n=1, s=1,
r=3$ and $n=2, s=3, r=3$.  In the first subcase $g^{-1}(p)$ splits
in two $A_3$-singularities, and in the second one, it splits in
four $D_4$-singularities.

The only remaining case is when $r\ge 4$ and
$\sum_{j=1}^r\varphi(\lambda_{i_j})\ne\varphi(\lambda_{i_k})$
whatever is $1\le k\le r$. Then $f^{-1}(E'_{p})$ splits into a
number of copies, denote one of them by $C$, of $2^m$-sheeted
Galois covering of $\mathbb P^1=E'_p$ branched at $r$ points,
where $m\geq 1$. By the Hurwitz formula,
$$
g(C)=2^{m-2}r-2^m+1\geq 1,
$$
which implies that the singular points $q\in g^{-1}(p)$ are not
canonical.\qed

\begin{lem}\label{k3}
Suppose that the Galois group of the covering $Y\to\mathbb P^2$
is $(\mathbb Z /2\mathbb Z)^3$ and that the
line arrangement $\mathcal L=L_0\cup\dots\cup L_6$ have no
$r$-fold singular points with $r\geq 4$. If $\varphi(\lambda_i)\ne
0$ for any $0\le i\le 6$ and there are two distinct lines
$L_{i_1}$ and $L_{i_2}$ with $\varphi
(\lambda_{i_1})=\varphi(\lambda_{i_2})$, then $p_g(X)\neq 0$.
\end{lem}

\proof By (\ref{new}), $Y$ can be given by equations
$$
u_j^2=\prod l_i(v_1,v_2)^{a_{i,j}}, \qquad j=1,2,3,
$$
where $(a _{i,1},a_{i,2} ,a _{i,3})=\varphi (\lambda _i)$.
Since $\varphi$ is an epimorphism to $(\mathbb Z /2\mathbb Z)^3$,
there are at most four lines with equal values of $\varphi$, so
that up to renumbering of lines and acting on $\varphi$ by
an automorphism of
$(\mathbb Z /2\mathbb Z)^3$
there are four cases to consider:
\begin{itemize}
\item (four equal values) $\varphi(\lambda_1)=\varphi(\lambda_2)
=\varphi(\lambda_3)=\varphi(\lambda_4)=(1,0,0)$,
$\varphi(\lambda_5)=(0,1,0)$, and $\varphi(\lambda_6)=(0,0,1)$;
\item (three equal values)
$\varphi(\lambda_1)=\varphi(\lambda_2)=\varphi(\lambda_3)=(1,0,0)$,
$\varphi(\lambda_4)=(0,1,0)$, and $\varphi(\lambda_5)=(0,0,1)$.
\item (two pairs of equal values)
$\varphi(\lambda_1)=\varphi(\lambda_2)=(1,0,0)$,
$\varphi(\lambda_3)=\varphi(\lambda_4)=(0,1,0)$, and
$\varphi(\lambda_5)=(0,0,1)$. \item (one pair of equal values)
$\varphi(\lambda_1)=\varphi(\lambda_2)=(1,0,0)$,
$\varphi(\lambda_3)=(0,1,0)$, $\varphi(\lambda_4)=(0,0,1)$, while
$\varphi(\lambda_i)$ with $i\in\{ 0,5,6\}$ are distinct from each
other and distinct from $(1,0,0),(0,1,0),(0,0,1)$.

\end{itemize}

In the first three cases the function $u=u_1u_2u_3\in \mathbb
C(Y)$ satisfies the following equation
\begin{equation} \label{ccc} u^2=l_1(v_1,v_2)\dots l_5(v_1,v_2)l_6(v_1,v_2)^a,
\end{equation}
where $a=0$ or $1$ (in the first case, $a=1$). Such an equation
defines a double covering $Z\to \mathbb P^2$ branched in six lines
($L_1,\dots, L_6$ if $a=1$ and $L_1,\dots, L_5, L_0$ if $a=0$).
Since the line arrangement has no $r$-fold points with $r\geq 4$,
$Z$ has only canonical singularities, and therefore it is a
$K3$-surface, which implies $p_g(Z)=1$. The inequality $p_g(X)\geq
1$ follows now from the existence of a dominant rational map from
$X$ to $Z$.

To
complete the proof, let us notice that the fourth case is
impossible. Indeed, it is impossible to satisfy the relation
$\varphi(\lambda_0)+\varphi(\lambda_5)+\varphi(\lambda_6)=(0,1,1)$,
by three distinct elements among $(1,1,0)$, $(0,1,1)$, $(1,0,1)$,
and $(1,1,1)$.\qed

\subsection{Campedelli surfaces as Galois coverings branched over Campedelli
arrangements}\label{asbranched} Let $\mathcal L$ be a line arrangement in $\mathbb
P^2$ consisting of seven distinct lines $L_\alpha$ labeled by the
non-zero elements $\alpha\in (\mathbb Z/2\mathbb Z)^3$. We call
such a labeled arrangement $\mathcal L$ a {\it Campedelli line
arrangement} if it has neither $r$-fold points with $r\geq 4$ nor
triple points $p_{\alpha_1,\alpha_2,\alpha_3}= L_{\alpha_1}\cap
L_{\alpha_2}\cap L_{\alpha_3}$  with
$\alpha_1+\alpha_2+\alpha_3=0$. We say that a Campedelli line
arrangement $\mathcal L=\sum L_{\alpha}$ is obtained from a
Campedelli line arrangement $\mathcal L'=\sum L_{\alpha}'$ by
means of {\it renum\-ber\-ing}  of lines if there is an
automorphism $\tau \in \text{Aut}(\mathbb Z/2\mathbb Z)^3$ such
that $L_{\alpha}=L_{\tau(\alpha)}'$ for any $\alpha \in (\mathbb
Z/2\mathbb Z)^3\setminus \{ 0\}$.

Given a Campedelli line arrangement $\mathcal L$, one can consider
the Galois covering $Y(\mathcal L)\to \mathbb P^2$ with Galois
group $(\mathbb Z/2\mathbb Z)^3$ branched in $\mathcal L$ and
defined by the epimorphism $\varphi :H_1(\mathbb P ^2 \setminus
\mathcal L, \mathbb Z )\to  (\mathbb Z/2\mathbb Z)^3$ given by
$\varphi(\lambda_{\alpha})=\alpha$. We call this covering {\it the
Galois covering  branched over a Campedelli arrangement} $\mathcal
L$. Clearly, a renumbering of a Campedelli arrangement leads to an
equivalent covering.

\begin{thm}\label{thm1} For any
Campedelli surface $X$ there exists a Campedelli line arrangement
$\mathcal L$ such that $X_ {can}= Y(\mathcal L)$.
\end{thm}

\proof By Corollary \ref{claim}, given a Campedelli surface $X$
there exists an arrangement $\mathcal L$ of seven distinct lines
in $\mathbb P^2$ such that $X_{can}$ is a $(\mathbb Z/2\mathbb
Z)^3$-Galois covering of $\mathbb P^2$ branched in $\mathcal L$.
Since $X_{can}$ has only canonical singularities, Lemma
\ref{A-D-E} implies that $\mathcal L$ have no neither any $r$-fold
point with $r\geq 4$ nor any 3-fold point which is not a branch
point. Now Lemma \ref{k3} applies and shows that $\mathcal L$ is a
Campedelli arrangement. \qed

The following, converse, statement is proved in \cite{Ku}.

\begin{thm} \label{thm2} {\rm (\cite{Ku})}
For any Campedelli line arrangement $\mathcal L$ the surface
$Y(\mathcal L)$ is isomorphic to the canonical model of a
Campedelli surface. \qed \end{thm}

If a Campedelli line arrangement $\mathcal L$ has no triple
points, then by Lemma \ref{lem 1.1}, the surface $Y(\mathcal L)$
is nonsingular (so that it is itself a Campe\-de\-l\-li surface,
$X=X_{can}$) and it can be imbedded as a complete intersection
into the weighted projective space
$$
\mathbb P^9_w=\mathbb P^9(1,1,1,2,2,2,2,2,2,2)
$$
with three weight-1 coordinates $z_i$, $i=0,1,2$, and seven
weight-2 coordinates $u_{\alpha}$, $\alpha\in (\mathbb Z/2\mathbb
Z)^3\setminus \{ 0\}$. Namely, in accordance with what was seen in
subsection \ref{linecovering},
$Y(\mathcal L)$ is isomorphic to a surface in $\mathbb P^9_w$ given
by
\begin{equation} \label{equations}
\begin{array}{l}
u_{(1,0,0)}^2=l_{(1,0,0)}
l_{(1,1,0)}l_{(1,0,1)} l_{(1,1,1)} \\
u_{(0,1,0)}^2= l_{(0,1,0)}
l_{(1,1,0)} l_{(0,1,1)} l_{(1,1,1)} \\
u_{(0,0,1)}^2= l_{(0,0,1)}
l_{(0,1,1)} l_{(1,0,1)} l_{(1,1,1)} \\
u_{(1,1,0)}^2= l_{(1,0,0)}
l_{(0,1,0)} l_{(1,0,1)}l_{(0,1,1)} \\
u_{(1,0,1)}^2= l_{(1,0,0)}
l_{(0,0,1)} l_{(1,1,0)} l_{(0,1,1)} \\
u_{(0,1,1)}^2= l_{(0,1,0)}
l_{(0,0,1)} l_{(1,0,1)} l_{(1,1,0)} \\
u_{(1,1,1)}^2= l_{(1,0,0)}
l_{(0,1,0)} l_{(0,0,1)} l_{(1,1,1)} .
\end{array}
\end{equation}
where
$
l_{\alpha}(z_0,z_1,z_2)=0$ are linear equations of
$L_{\alpha}\subset \mathcal L$ in $\mathbb P^2$.

Note  that $u_{\alpha}$ satisfy the following relations
\begin{equation} \label{dependness}
\begin{array}{l}
\displaystyle u_{(1,1,0)}=\frac{u_{(1,0,0)}u_{(0,1,0)}}{
l_{(1,1,0)} l_{(1,1,1)}},\, \, \,
u_{(1,0,1)}=\frac{u_{(1,0,0)}u_{(0,0,1)}}{
l_{(1,0,1)} l_{(1,1,1)}},\\
\\ \displaystyle
u_{(0,1,1)}=\frac{u_{(0,1,0)}u_{(0,0,1)}}{
l_{(0,1,1)} l_{(1,1,1)}}, \,\, \, \displaystyle
u_{(1,1,1)}=\frac{u_{(1,0,0)}u_{(0,1,0)}u_{(0,0,1)}}{
l_{(1,1,0)} l_{(1,0,1)} l_{(0,1,1)}
l_{(1,1,1)}}.
\end{array}
\end{equation}
Note also that if $\mathcal L'$ is obtained from $\mathcal L$ by a
renumbering of the lines $\mathcal L$ given by an automorphism
$\tau \in\text{Aut}(\mathbb Z/2\mathbb Z)^3$, then this
renumbering (in order to save the form of the equations in
(\ref{equations})) defines the renumbering of $u_{\alpha}$ by the
automorphism $\tau^{-1}$.

\subsection{Moduli space of the Campedelli surfaces} \label{ssection}
In this section, we identify the moduli space of Campedelli surfaces
with the moduli space of Campedelli line arrangements. Here and
further, we apply to Campedelli surfaces the following general
property of minimal surfaces of general type: their isomorphisms
(respectively, automorphisms) are in a natural bijection with the
isomorphisms (respectively, automorhisms) of their canonical models.

 As above, let a Galois covering $g: Y(\mathcal L) \to \mathbb P
^2$ with Galois group $G\simeq (\mathbb Z /2\mathbb Z)^{3}$ be
branched along a Campedelli line arrangement $\mathcal L=\sum
L_{\alpha}$, where the sum is taken over all $\alpha\in G$, $\alpha
\neq 0$, and be determined by an epimorphism $\varphi : H_1(\mathbb
P ^2 \setminus \mathcal L, \mathbb Z )\to G$ such that $\varphi
(\lambda_{\alpha})=\alpha$. Denote by $X=X(\mathcal L)$ the minimal
nonsingular model of $Y(\mathcal L)$ constructed in subsection
\ref{Yresolution}. Since $\mathcal L$ has neither $r$-fold points
with $r\geq 4$ nor triple points $p_{\alpha_1,\alpha_2,\alpha_3}=
L_{\alpha_1}\cap L_{\alpha_2}\cap L_{\alpha_3}$  with
$\alpha_1+\alpha_2+\alpha_3=0$, this construction reduces to the
composition $\sigma : \hat{\mathbb P}^2\to \mathbb P^2$ of the
blow-ups with centers at all the $3$-fold points of $\mathcal L$
followed by the covering $f:X(\mathcal L)\to \hat{\mathbb P}^2$
induced by the lift $\hat\varphi$ of $\varphi$.

Denote by $f_\sigma$ the composition $f_\sigma=\sigma\circ f :
X(\mathcal L)\to \mathbb P^2$.

\begin{lem} \label{can} {\rm (\cite{Ku})}
The bicanonical system $\vert 2K_X\vert $ of $X=X(\mathcal L)$ is
equal to $| f^*_\sigma L|$, where $L\subset \mathbb P^2$ is a line
in $\mathbb P^2$. \qed
\end{lem}

The next Lemma is a straightforward corollary of Proposition
\ref{0.1}.

\begin{lem} \label{top} Let $\mathcal L_1=\sum_{i=1}^7 L_{1,\alpha_i}$ and
$\mathcal L_2=\sum_{i=1}^7 L_{2,\beta_i}$, $\alpha_i,\beta_i\in
G=(\mathbb Z/2\mathbb Z)^3$, $\alpha_i,\beta_i \neq 0$, be two
Campedelli line arrangements in $\mathbb P^2$ such that
$L_{1,\alpha_i}=L_{2,\beta_i}$ for $i=1,\dots, 7$. Then the Galois
coverings $Y(\mathcal{L}_1)\to\mathbb P^2$ and
$Y(\mathcal{L}_2)\to \mathbb P^2$ are equivalent if, and only if,
$\mathcal L_1$ can be obtained from $\mathcal L_2$ by means of
renumbering of lines. \qed
\end{lem}

\begin{thm} \label{autcam} Let $X_{1,can}=Y(\mathcal L_1)$
and $X_{2,can}=Y(\mathcal L_2)$ be two Galois coverings
$g_i:X_{i,can}\to \mathbb P^2$ branched over Campedelli line
arrangements $\mathcal L_1$ and $\mathcal L_2$. If $X_{1,can}$ and
$X_{2,can}$ are isomorphic, then any isomorphism $\nu :X_{1,can}\to
X_{2,can}$ can be included in a commutative diagram

\begin{picture}(300,90)
\put(125,70){$X_{1,can}$} \put(165,77){$\nu $}
\put(155,75){\vector(1,0){40}} \put(202,70){$X_{2,can}$}
\put(131,65){\vector(0,-1){30}}\put(134,47){$g_1$}
\put(205,65){\vector(0,-1){30}}\put(210,47){$g_2$}
\put(128,20){${\mathbb P ^2}$} \put(170,13){$\psi$}
\put(145,25){\vector(1,0){50}} \put(201,20){$\mathbb P ^2$\,.}
\end{picture}
\end{thm}

\proof  Consider the resolutions $X_i=X(\mathcal L_i)$ of
$X_{i,can}=Y(\mathcal L_i)$, the associated coverings $f_i: X_i\to
\hat{\mathbb P}^2$, and the composed morphisms
$f_{\sigma,i}=\sigma\circ f_i: X_i\to\mathbb
P^2$. As it was mentioned above,
since $X_i$ are minimal surfaces of general type, any isomorphism
between their canonical models, $X_{1,can}\to X_{2,can}$ lifts
uniquely to an isomorphism $X_1\to X_2$, and vice versa. Thus, it is
sufficient to pick an isomorphism $\nu :X_1\to X_2$ and to find a
projective transformation $\psi$ such that $\psi\circ
f_{\sigma,1}=f_{\sigma,2}\circ\nu$. Moreover, the latter relation
would follow from the corresponding relation between the induced
maps of the function fields: $\nu^*\circ f^*_{\sigma,2}=
f^*_{\sigma,1}\circ\psi^*$.

As for any Campedelli surface, the torsion subgroup $\mbox{\rm
Tors}(X_i)$ of $H^2(X_i,\mathbb Z)$ is $2$-torsion and isomorphic to
$(\mathbb Z /2\mathbb Z)^{3}$. Given $\alpha\in \mbox{\rm Tors}\,
(X_i)$, $\alpha\neq 0$, the linear system $|K_{X_i}+\alpha |$ is
non-empty as it follows from Serre duality,
$$
\dim H^2(X_i, \mathcal O_{X_i}(K_{X_i}+\alpha))=\dim H^0(X_i,
\mathcal O_{X_i}(\alpha))=0,
$$
and the Riemann-Roch theorem. Hence, there exists at least one
effective divisor $D_{\alpha}\in |K_{X_i}+\alpha |$, and
$2D_{\alpha}\in |2K_{X_i}|$. Since $X_i$ are minimal surfaces of
general type, we have $ \dim H^0(X_i,\mathcal
O(2K_{X_i}))=K^2_{X_i}+1=3$. On the other hand, $ \dim H^0(\mathbb
P^2,\mathcal{O}_{\mathbb P^2}(L))=3$, where $L$ is a line in
$\mathbb P^2$, while by Lemma \ref{can} we have $\vert 2K_{X_i}\vert
= \vert f^*_\sigma (L)\vert $. Finally, $|2K_{X_i}|
=f_{\sigma,i}^*(|L|)$ and $D\in |K_{X_i}+\alpha |$ for some $\alpha
\in \mbox{\rm Tors}(X_i)$ if and only if $2D=
f_{\sigma,i}^*(\widetilde L)$ for some $\widetilde L \in |L|$.

The only lines $\widetilde L \in |L|$ for which the divisors
$f_{\sigma,i}^*(\widetilde L)$ are divisible by $2$ are the seven
branch lines belonging to $\mathcal L_i$. Hence, they give all the
different torsion elements and can be relabeled by the torsion
elements so that $\mathcal L_i=\sum L_{i,\alpha}$, where the sum is
taken over the nonzero torsion elements, and
$\frac{1}{2}f_{\sigma,i}^*(L_{i,\alpha})=D_{i,\alpha}\in
|K_{X_i}+\alpha |$. (Note that this labeling of lines may not
coincide with the initial one.)

Let $\nu:X_1 \to X_2$ be an isomorphism. It induces an isomorphism
of torsion groups, $\nu^*: \mbox{\rm Tors}\, (X_2)\to \mbox{\rm
Tors}\, (X_1)$, and isomorphisms of linear systems,
$$\nu^*: H^0(X_2, \mathcal
O_{X_2}(K_{X_2}+\alpha)) \to H^0(X_1, \mathcal
O_{X_1}(K_{X_1}+\nu^*(\alpha)))$$ for each $\alpha\in \mbox{\rm
Tors}\, (X_2)$. Therefore, $\nu^*(D_{2,\alpha})=D_{1,\nu^*(\alpha
)}$ for any $\alpha \in \mbox{\rm Tors}(X_2)$, $\alpha \neq 0$, and
we get
$$
\begin{array}{l}
\nu^*(f_{\sigma,2}^*(L_{2,\alpha_1}-L_{2,\alpha_2}))=\nu^*(2D_{2,\alpha_1}-2D_{2,\alpha_2})=
2D_{1,\nu^*(\alpha_1)}-2D_{1,\nu^*(\alpha_2)}=
\\ \hspace{3.8cm} f_{\sigma,1}^*(L_{1,\nu^*(\alpha_1)}-L_{1,\nu^*(\alpha_2)})
\end{array}
$$
for any non zero $\alpha_1,\alpha_2\in \mbox{\rm Tors}\, (X_2)$.
Since any rational function is defined uniquely up to multiplication
by a constant by its divisors of zeros and poles, it implies the
existence of a system of constants $c_{\alpha_1,\alpha_2}$ such that
\begin{equation}
\label{function}\displaystyle
\nu^*(f_{\sigma,2}^*(\frac{l_{2,\alpha_1}(v_{1},v_{2})}
{l_{2,\alpha_2}(v_{1},v_{2})}))=
c_{\alpha_1,\alpha_2}f_{\sigma,1}^*(\frac{l_{1,\nu^*(\alpha_1)}(v_{1},v_{2})}
{l_{1,\nu^*(\alpha_2)}(v_{1},v_{2})}),
\end{equation}
where $v_1,v_2$ are affine coordinates in $\mathbb P^2$ and
$l_{2,\alpha}, l_{1,\beta}$ are linear equations of the
corresponding lines. Since the functions
$f_{\sigma,i}^*(\frac{l_{i,\alpha_1}(v_{1},v_{2})}
{l_{i,\alpha_2}(v_{1},v_{2})}))$ generate the subfields
$f_{\sigma,i}^*(\mathbb C(\mathbb P^2))$ of $\mathbb C(X_i)$, the
relations (\ref{function}) imply the existence of a projective
transformation $\psi:\mathbb P^2\to\mathbb P^2$ such that
$f_{\sigma,1}^*\circ \psi^*=\nu^*\circ f_{\sigma,2}^*$.
 \qed

\begin{cor} \label{autcam+}
If $X=X(\mathcal L)$, where $\mathcal L$ is a generic Campedelli
line arrangement, then $\text{Aut}(X)=\text{Gal}(Y(\mathcal L)\to
\mathbb P^2)\simeq(\mathbb Z /2\mathbb Z)^{3}$. \qed
\end{cor}

Denote by $\mathcal P=\mathbb P^2\times \dots \times \mathbb P^2$
the product of seven copies of the projective plane.  We consider
each factor in this product as the dual projective plane, so that
elements of each factor are lines in the initial $\mathbb P^2$. In
addition, we numerate the factors of $\mathcal P$ by the non-zero
elements $\alpha \in G=(\mathbb Z/2\mathbb Z)^3$. Let $D$ be the
union of all diagonals in $\mathcal P$,
$$ \begin{array}{ll} T_3=  \{ &
\mathcal L\in \mathcal P\, \mid \,
\exists\, \alpha_{i_1},\alpha_{i_2},\alpha_{i_3}\,\, \text{such
that} \\ &  \alpha_{i_1}+\alpha_{i_2}+\alpha_{i_3}=0\, \,
\text{and} \, \, L_{\alpha_{i_1}}\cap L_{\alpha_{i_2}}\cap
L_{\alpha_{i_3}}\neq \emptyset \} , \end{array}
$$
$$ \begin{array}{ll} T_4=  \{ &
\mathcal L\in \mathcal P\, \mid \,
\exists\, \alpha_{i_1},\alpha_{i_2},\alpha_{i_3}, \alpha_{i_4}\,\,
\text{such that} \\ & L_{\alpha_{i_1}}\cap L_{\alpha_{i_2}}\cap
L_{\alpha_{i_3}}\cap L_{\alpha_{i_4}}\neq \emptyset \} .
\end{array}
$$
The group  $\text{PGL}(2,\mathbb C)\times \text{Aut}\, G$  acts on
$\mathcal P\setminus(D\cup T_3\cup T_4)$ as follows:
$\text{PGL}(2,\mathbb C)$ acts in a usual way on each factor of
$\mathcal P$, and the elements $h$ of $\text{Aut}\, G$ permute the
factors, $h:\mathbb P^2_{\alpha}\to \mathbb P^2_{h(\alpha)}$.

The following theorem is a consequence of Lemma \ref{top} and
Theorems  \ref{thm1}, \ref{thm2}, and \ref{autcam}.
\begin{thm}\label{modcamp}
The moduli space $\mathcal M$ of the Campedelli surfaces is
isomorphic to the quotient space $$(\mathcal P\setminus(D\cup
T_3\cup T_4))/\text{PGL}(2,\mathbb C)\times \text{Aut}\, G.$$ \qed
\end{thm}

Note that, as a result, all Campedelli surfaces are deformation
equivalent.

\section{Real Campedelli surfaces}

\subsection{An extension of the automorphism group}\label{Kl}
For any complex space $X$, denote by $\text{Kl}=\text{Kl}(X)$ the
group of holomorphic and anti-holomorphic bijections $X\to X$.
Recall that, by definition, an anti-holomoprhic map $X\to X$ can
be seen as a holomorphic map $X\to \overline X$, where $\overline
X$ states for the complex conjugate to $X$.

Note ({\it cf.} subsection \ref{ssection}) that for any minimal
surface $X$ of general type the groups $\text{Kl}(X)$ and
$\text{Kl}(X_{can})$ are naturally isomorphic. In what follows we
identify them as soon as it does not lead to a confusion.

Clearly, if $\text{Kl}$ contains at least one anti-holomorphic
element, the holomorphic elements form in $\text{Kl}$ a subgroup
$\text{Aut}=\text{Aut}(X)$ of index $2$. In other words, there is
a short exact sequence $1\to\text{Aut}\to\text{Kl}\to H\to 1$,
where $H\simeq \mathbb Z/2$ or $1$. We denote by
$\text{kl}:\text{Kl}\to H$ the homomorphism of this sequence.

The real structures on $X$ are the elements $c\in\text{Kl}(X)$
such that $\text{kl}(c)\neq 1$ and $c^2=\text{id}$. Two real
structures, $c_1$ and $c_2$ are called {\it equivalent} (or {\it
isomorphic}) if there exists $h\in\text{Aut}(X)$ such that $h\circ
c_2= c_1\circ h$. Recall that on the projective plane $\mathbb
P^2$ (as well as on any projective space of even dimension) any
two real structures are equivalent by a projective transformation.

\subsection{A criteria of existence of real structures on Campedelli surfaces}
Given a Galois covering $g:Y(\mathcal L)\to {\mathbb P}^2$ branched
over a Campedelli line arrangement $\mathcal L$ and the associated
Campedelli surface $X=X(\mathcal L)$ together with the composed map
$f_\sigma=\sigma\circ f : X \to\mathbb P^2$, we say that $c_X\in
\text{Kl}(X)$ is {\it lifted} from ${\mathbb P}^2$ if there exists
$c_P\in \text{Kl}({\mathbb P}^2)$ such that the following diagram is
commutative

\begin{picture}(300,90)
\put(127,70){$X$} \put(165,78){$c_X$} \put(141,75){\vector(1,0){50}}
\put(200,70){$X$}
\put(131,65){\vector(0,-1){30}}\put(134,47){$f_{\sigma}$}
\put(205,65){\vector(0,-1){30}}\put(210,47){$f_{\sigma}$}
\put(128,20){${\mathbb P}^2$} \put(170,17){$c_P$}
\put(145,25){\vector(1,0){50}} \put(201,20){${\mathbb P}^2.$}
\end{picture}

\begin{thm} \label{realthm}
For any Campedelli line arrangement $\mathcal L$, every
$c_X\in\text{Kl} (X)$ is lifted from ${\mathbb P}^2$. In particular,
if $X$ has a real structure $c_X$, then there exists a real
structure $c_P$ on ${\mathbb P}^2$ such that $c_P\circ
f_{\sigma}=f_{\sigma}\circ c_X$.
\end{thm}

\proof If $c_X\in \text{Aut}(X)$, then $c_X$ is lifted from
${\mathbb P}^2$ by Theorem \ref{autcam}. Let $c_X\in \text{Kl}(X)$
and $c_X\not\in \text{Aut}(X)$. Then $c_X:X\to\overline X$ is a
holomorphic isomorphism. Consider the complex conjugated covering
$\overline f_{\sigma}: \overline X\to \overline{{\mathbb P}^2}$. By
Theorem \ref{autcam}, there is a holomorphic isomorphism $c_P:
{\mathbb P}^2\to \overline{{\mathbb P}^2}$ which makes commutative
the following diagram

\begin{picture}(300,90)
\put(125,70){$X$} \put(165,79){$c_X$} \put(145,75){\vector(1,0){50}}
\put(202,70){$\overline X$}
\put(131,65){\vector(0,-1){30}}\put(134,47){$f_{\sigma}$}
\put(205,65){\vector(0,-1){30}}\put(210,47){$\overline f_{\sigma}$}
\put(126,20){ ${\mathbb P}^2$} \put(160,17){ ${c_P}$}
\put(145,25){\vector(1,0){50}} \put(201,20){$\overline{
{\mathbb P}^2}.$}
\end{picture}

To get the last statement, it is sufficient to notice that
$c_P^2=\text{id} $ if $c_X^2=\text{id} $. \qed

\begin{cor}\label{realcor}
For any Campedelli line arrangement $\mathcal L\subset \mathbb
P^2$, the Campedelli surface $X=X(\mathcal L)$ admits a real
structure if, and only if, for a suitably  chosen real structure
$c_P$ of $\mathbb P^2$ the (labeled) Campedelli line arrangement
$\mathcal L$ is real, that is, there exists an automorphism
(renumbering) $\tau:(\mathbb Z/2\mathbb Z)^3\to(\mathbb Z/2\mathbb
Z)^3$ such that $c_P(L_\alpha)=L_{\tau(\alpha)}$ for each
$\alpha\in(\mathbb Z/2\mathbb Z)^3$, $\alpha\ne 0$.
\end{cor}

\proof In the case of a real arrangement, to lift $c_P$ it is
sufficient to notice that $c_P$ (as any real structure on $\mathbb
P^2$) has a whole real projective plane of fixed points, to pick
such a fixed point in the complement of the arrangement, and to
identify the unbranched points of $X_{can}$ with classes of pathes
issued from the fixed point so that $c$ becomes properly acting on
$X_{can}$. A renumbering induced by a transformation of $\mathbb
P^2$ is a homomorphism, since it factors through the induced action
on $H_1(\mathbb P ^2 \setminus \mathcal L, \mathbb Z )$. \qed

\subsection{Real Campedelli line arrangements}
The Galois group $G=\text{Gal}(X/ \hat{\mathbb P}^2)\simeq (\mathbb
Z/2\mathbb Z)^3$ is a subgroup of $\text{Aut}(X)$.  As it follows
from Theorem \ref{realthm}, $G$ is a normal subgroup of
$\text{Kl}(X)$, and in addition, by Corollary \ref{realcor},
$c(L_{\alpha} )=L_{c\alpha c^{-1}}$ for any $\alpha \in G$ and
$c\in\text{Kl}(X)$.

\begin{prop}\label{raltype}
Let $\mathcal L$ be a Campedelli line arrangement which is real with
respect to some real structure $c_P:\mathbb P^2\to\mathbb P^2$. Then
either $\mathcal L$ consists of seven real lines or it consists of
three real lines and two pairs of complex conjugated lines.
Respectively, $c_P$ acts on the labeling of $\mathcal L$ either
identically or not.
\end{prop}

\proof The homomorphism
$\alpha \in G=(\mathbb Z/2\mathbb Z)^3\mapsto c\alpha
c^{-1}\in G=(\mathbb Z/2\mathbb Z)^3$, where $c$ is the
real structure on $X$, is an involution, and, as any involution on a
$\mathbb Z/2$-vector space, it splits into irreducible $1$- and
$2$-dimensional components. In dimension $3$, there are only two
possibilities, either the involution is trivial or it contains a
$2$-dimensional irreducible component, that is an involution
interchanging two generators. In the first case, all
$\alpha$ are fixed, and hence all the lines are real. In the second
case, there are three and only three fixed elements, and hence three
and only three real lines. \qed

Let call a Campedelli line arrangement  $\mathcal L$ {\it purely
real}  if it consists of seven real lines and {\it mixed real} if
it consists of three real lines and two pairs of complex
conjugated lines.

Given a real structure $c_X$, denote by $\text{Kl}(X,c_X)$ the
subgroup of $\text{Kl}(X)$ generated by $G$ and $c_X$. If
$X=X(\mathcal L)$ and $L$ is real with respect to a real structure
$c_P$ on $\mathbb P^2$, then the subgroup $\text{Kl}(X,c_X)$ does
not depend on the choice of a lift $c_X$ of $c_P$ and we denote it
by $\text{Kl}(X,c_P)$. Note that for a generic real Campedelli
line arrangement $\mathcal L$ it holds $\text{Aut} X(\mathcal
L)=G$, so that $\text{Kl(X)}=\text{Kl}(X,c_X)$ for any $c_X$.

\begin{prop}
\label{klcam} Let $X=X(\mathcal L)$ be a Campedelli surface
associated with a Campedelli line arrangement $\mathcal L$ which
is real with respect to $c_P$. Then:
\begin{itemize}
\item[($i$)] if $\mathcal L$ is a purely real line arrangement,
then $\text{Kl}(X, c_P)\simeq (\mathbb Z/2\mathbb Z)^4$;
and if $\mathcal L$ is a generic purely real line arrangement, then there are
exactly eight different real structures on $X$;
\item[($ii$)] if
$\mathcal L$ is a mixed real line arrangement, then
$\text{Kl}(X, c_P)\simeq \mathbb H\times(\mathbb Z/2\mathbb Z)$, where
$\mathbb H$ is the quaternion group of order eight; and
if $\mathcal L$ is a generic mixed real line arrangement, then there are
exactly four different real structures on $X$.
\end{itemize}
\end{prop}
\proof Pick a real point $p\in\mathbb P^2\setminus\mathcal L$ and
consider a real structure $c\in\text{Kl}(X) $ which lifts $c_P$
from $\mathbb P^2$ to $X$ and have fixed points over $p$. If all
the lines are real, then $c \alpha c^{-1}=\alpha $ for any $\alpha
\in G$ (indeed, since $c=\text{id}$ at each point of the $G$-orbit
over $p$, the relation $c \alpha c^{-1}=\alpha $ holds at the
points of this $G$-orbit, and, hence, it holds everywhere). If
there are only three real lines in the arrangement, then in a
suitable basis $e_1,e_2,e_3$ of $G$ the (renumbering) involution
$\alpha \mapsto c \alpha c^{-1}$ acts as $e_1\mapsto e_2$ and
$e_3\mapsto e_3$. Therefore, in the latter case, $\text{Kl}(X,
c_P)$ splits in a direct sum of $\mathbb Z/2$ generated by $e_3$
with a non-commutative group of order 8 generated by $e_1, e_2$,
and $c$.

Since for a generic arrangement it holds $\text{Kl(X)}=\text{Kl}(X,c_X)$,
the statements concerning the generic cases follow now from
enumerating anti-involutions in $\text{Kl}(X, c_P)\simeq (\mathbb Z/2\mathbb Z)^4$
and, respectively, $\text{Kl}(X, c_P)\simeq \mathbb H\times(\mathbb Z/2\mathbb Z)$. \qed

\subsection{Purely real Campedelli line arrangements}\label{purely}  Let
$\mathcal L$ be a Cam\-pe\-del\-li line arrangement $\mathcal L=
\cup L_{\alpha}$ which is purely real with respect to a real
structure $c_P: \mathbb P^2\to\mathbb P^2$. Choose homogeneous
coordinates $(z_0,z_1,z_2)$ in $\mathbb P^2$ such that $c_P$ turns
in the standard complex conjugation
$$
c_P(z_0,z_1,z_2)=(\overline z_0,\overline z_1,\overline z_2).$$ Then
each of the lines $L_{\alpha}\in\mathcal L$, $\alpha\in G\setminus
\{ 0\}$, is given by equation
$$a_{\alpha,0}z_0+a_{\alpha,1}z_1+a_{\alpha,2}z_2=0$$
with real coefficients, $a_{\alpha,i}\in \mathbb R$.

Consider the set $\mathbb R\mathbb P^2=\{ (z_0,z_1,z_2) \mid \,
z_i\in \mathbb R\}$ of real points of $\mathbb P^2$.  If $\mathcal
L$ has no triple points, then $\mathcal L$ divides $\mathbb
R\mathbb P^2$ into twenty two $n$-gons $P_i$, $i=1,\dots,22$,
$3\leq n\leq 7$. The collection $(m_3,\dots,m_7)$, where $m_n$ is
the number of $n$-gons $P_i$, is called the {\it type} of
$\mathcal L$.

The following description of topology of the inverse image of $P_i$
in the associated Campedelli surface $X(\mathcal L)$ is a
straightforward consequence of the construction of ramified
coverings.

\begin{prop}\label{gluing} For any polygon $P_i$ of a purely
real Campedelli line arrangement $\mathcal L$ without triple
points, its inverse image $f^{-1}(P_i)\subset X(\mathcal L)$ is a
two-manifold and it is homeomorphic to the following quotient of
$P_i\times G$, $G=(\mathbb Z/2\mathbb Z)^3$: the points
$(a,\beta)$ and $(b,\gamma)$ are identified if $a=b\in L_\alpha$
where $\gamma=\beta+\alpha$.\qed
\end{prop}
A triangle $P_i$ bounded by $L_{\alpha_1}$, $L_{\alpha_2}$, and
$L_{\alpha_3}$ is said to have {\it linear {\rm(}in{\rm)}dependent
sides}, if $\alpha_1$, $\alpha_2$, $\alpha_3$ are linear
(in)dependent.

\begin{cor}\label{lll}
For any $n$-gon $P_i$ of a purely
real Campedelli line arrangement $\mathcal L$ without triple
points,
\begin{itemize}
\item[($i$)] the Euler characteristic of $f^{-1}(P_i)$  is equal
to $8-2n$;
\item[($ii$)]
$f^{-1}(P_i)$ is the disjoint union of two copies of
$\mathbb R\mathbb P^2$, if $n=3$ and the triangle $P_i$ has linear depended
sides;
\item[($iii$)] $f^{-1}(P_i)$ is the
two-dimensional sphere, if $n=3$ and the triangle $P_i$ has
linear independent sides;
\item[($iv$)] $f^{-1}(P_i)$ is connected, if $n=4$,
and it is
orientable  if, and only if, $\alpha_1+\dots +\alpha_4=0$, where $\alpha_j$ are
the labels of the sides $L_{\alpha_j}$ of $P_i$;
\item[($v$)]
$f^{-1}(P_i)$ is a
connected non-orientable two-manifold, if $n\geq 5$.
\end{itemize}
\end{cor}

\proof The Euler characteristic $e(f^{-1}(P_i))$ is equal to
$$e(f^{-1}(P_i))=8-4n+2n=8-2n$$
according to the cellular decomposition given by Proposition
\ref{gluing}.

Let $L_{\alpha_1},\dots, L_{\alpha_n}$ be the sides of $P_i$.
Consider a subgroup $G_{P_i}=\langle
\alpha_1,\dots,\alpha_n\rangle$ of $G$ generated by
$\alpha_1,\dots,\alpha_n$. As it follows from
Proposition \ref{gluing}, the number of connected components of
$f^{-1}(P_i)$ coincides with the index of $G_{P_i}$ in $G$. On the
other hand, since $n>2$, either $G_{P_i}$ coincides with $G$ or it
is a subgroup of index $2$, and in the latter case, $P_i$ is a
triangle with linear dependent sides. Therefore, $f^{-1}(P_i)$ is
connected except in the case of triangles with linear dependent
sides and, moreover, if $P_i$ is a triangle with linear dependent
sides, then $f^{-1}(P_i)$ consists of two connected components.

If $n=3$, then $e(f^{-1}(P_i))=2$. Hence, if $P_i$ is a
triangle with linear independent sides, then $f^{-1}(P_i)$ is the
$2$-sphere, and if $P_i$ is a triangle with linear dependent sides, then
$f^{-1}(P_i)$ is the disjoint union of two copies of $\mathbb
R\mathbb P^2$.

Let $n\geq 4$. Then, $P_i$ has three successive sides whose indices
$\alpha_1,\alpha_2,\alpha_3$ are linear independent. After
renumbering we can assume that $\alpha_1=(1,0,0)$,
$\alpha_2=(0,1,0)$, and $\alpha_3=(0,0,1)$. Following
Proposition \ref{gluing}, perform a partial gluing of eight copies
$P_{\beta}=(P_i,\beta)$ of $P_i$ as is depicted in Fig. 1 (in Fig.
1, we denote the union of sides $L_{\alpha_4}\cup\dots\cup
L_{\alpha_n}$ by $\widetilde{L}_{\alpha}$).

\begin{picture}(300,180)
\put(0,150){\vector(0,-1){75}} \put(0,0){\vector(0,1){75}}
\put(75,75){\vector(0,1){75}} \put(75,75){\vector(0,-1){75}}
\put(150,0){\vector(0,1){75}} \put(150,150){\vector(0,-1){75}}
\put(225,75){\vector(0,1){75}} \put(225,75){\vector(0,-1){75}}
\put(300,0){\vector(0,1){75}} \put(300,150){\vector(0,-1){75}}
\put(0,75){\vector(1,0){75}}\put(150,75){\vector(-1,0){75}}
\put(150,75){\vector(1,0){75}} \put(300,75){\vector(-1,0){75}}
\put(75,150){\vector(-1,0){75}} \put(75,150){\vector(1,0){75}}
\put(225,150){\vector(-1,0){75}} \put(225,150){\vector(1,0,){75}}
\put(75,0){\vector(-1,0){75}} \put(75,0){\vector(1,0){75}}
\put(225,0){\vector(-1,0){75}} \put(225,0){\vector(1,0,){75}}
\put(30,-13){$\widetilde{L}_{\alpha}$}\put(105,-13){$\widetilde{L}_{\alpha}$}
\put(185,-13){$\widetilde{L}_{\alpha}$}\put(255,-13){$\widetilde{L}_{\alpha}$}
\put(30,153){$\widetilde{L}_{\alpha}$}\put(105,153){$\widetilde{L}_{\alpha}$}
\put(185,153){$\widetilde{L}_{\alpha}$}\put(255,153){$\widetilde{L}_{\alpha}$}
\put(30,78){${L}_{(1,0,0)}$}\put(105,78){${L}_{(1,0,0)}$}
\put(185,78){${L}_{(1,0,0)}$}\put(255,78){${L}_{(1,0,0)}$}
\put(3,128){${L}_{(0,1,0)}$}\put(78,128){${L}_{(0,0,1)}$}
\put(153,128){${L}_{(0,1,0)}$}\put(228,128){${L}_{(0,0,1)}$}
\put(303,128){${L}_{(0,1,0)}$}
\put(3,48){${L}_{(0,1,0)}$}\put(78,48){${L}_{(0,0,1)}$}
\put(153,48){${L}_{(0,1,0)}$}\put(228,48){${L}_{(0,0,1)}$}
\put(303,48){${L}_{(0,1,0)}$}
\put(25,105){$P_{(1,1,1)}$}\put(25,30){$P_{(0,1,1)}$}
\put(100,105){$P_{(1,1,0)}$}\put(100,30){$P_{(0,1,0)}$}\put(180,105){$P_{(1,0,0)}$}
\put(180,30){$P_{(0,0,0)}$}\put(250,105){$P_{(1,0,1)}$}\put(250,30){$P_{(0,0,1)}$}
\end{picture}
\\

\begin{center}
Fig. 1
\end{center}

Let $n=4$. Then, for $\widetilde{L}_{\alpha}=L_{\alpha_4}$ there are
four cases: either $\alpha_4=(1,1,0)$, or $\alpha_4=(1,0,1)$, or
$\alpha_4=(0,1,1)$, or $\alpha_4=(1,1,1)$. It is easy to see from
Fig. 1 that $f^{-1}(P)$ is non-orientable in the first three cases
and it is orientable in the last case.

Let, finally, $n\geq 5$. Then,
$\widetilde{L}_{\alpha}=L_{\alpha_4}\cup \dots \cup L_{\alpha_n}$
and at least one of $\alpha_4,\dots,\alpha_n$, say $\alpha_j$, has
to be equal to either $(1,1,0)$, or $(1,0,1)$, or $(0,1,1)$.
Therefore, the gluing  of $P_{(0,0,0)}$ and $P_{\alpha_j}$ along
$L_{\alpha_j}$ gives rise to non-orientability of $f^{-1}(P)$. \qed

Consider a real structure $c_X: X(\mathcal L)\to X(\mathcal L)$
which is a lift of $c_P$. According to Proposition \ref{klcam},
$c_X$ commutes with every element of $G$. Therefore, for any
$P_i$, $1\le i\le 22$, there exists one and only one $g_i\in G$
such that $c_X(x)=g_i(x)$ for any $x\in X$ with $f_\sigma(x)\in
P_i$. Using the same identification of $G$ with $(\mathbb
Z/2\mathbb Z)^3$ which we have already fixed introducing the
labeling of $\mathcal L$, $\varphi : H_1(\mathbb P ^2 \setminus
\mathcal L, \mathbb Z )\to (\mathbb Z/2\mathbb Z)^3$, we put
$$
g_i=({g}_{i,1},
{g}_{i,2},
{g}_{i,3})
$$
and introduce sign-triples
$$
Sign(P_i)=Sign_i=(
{sign}_{i,1},
{sign}_{i,2},
{sign}_{i,3}),
$$
where, by definition, ${sign}_{i,k}=(-1)^{g_{i,k}}, 1\le k\le 3$.
When we renumber the lines in $\mathcal L$ by means of an
automorphism $h:(\mathbb Z/2\mathbb Z)^3\to (\mathbb Z/2\mathbb Z)^3$
the labels $g_i$ of $P_i$ transform in $h(g_i)$; in particular,
the labels $Sign_i=(+,+,+)$ (corresponding to $g_i=0$) remain
unchanged under any renumbering. We call
{\it positive} the polygons $P_i$ with labels $Sign_i=(+,+,+)$.

The labels $Sign_i$ satisfy the following transition rule:
\begin{equation}\label{transition}
{sign}_{i,k}=(-1)^{a_{k}}sign_{j,k}
\end{equation}
if $P_i$ and $P_j$ have a common side on $L_\alpha$, $\alpha=(a_1,a_2,a_3)$.
In particular, if one of $Sign_i$ is given, then it determines all the other.

Let us notice that we switch from $g_i$ to $Sign_i$ by two reasons: first, it allows
us to distinguish more easily (say, on Figures) a labeling of lines,
$L_\alpha\mapsto \alpha$, from a labeling of polygons, $P_i\mapsto Sign_i$;
second, these signs have a natural meaning described
below (and are convenient in use).

Namely, to give an equivalent description of the above
sign-labeling, let consider the embedding of $Y(\mathcal L)$ into
$\mathbb P^9_{w}$ given by equations (\ref{equations}) and the
products
\begin{equation}\label{basicsigns}
\begin{array}{l}
l_{(1,0,0)}l_{(1,1,0)}l_{(1,0,1)}l_{(1,1,1)},\\
l_{(0,1,0)}l_{(1,1,0)}l_{(0,1,1)}l_{(1,1,1)},\\
l_{(0,0,1)}l_{(1,0,1)}l_{(0,1,1)}l_{(1,1,1)} \\
\end{array}
\end{equation}
participating in the first three equations (see subsection
\ref{asbranched} for notations related with $\mathbb P^9_{w}$). As
any homogeneous form of even degree with real coefficients, each
of the products has a well defined sign at any point of $\mathbb
RP^2$, where the product is nonzero. In particular, all the three
products have a well defined sign at the interior of each of $P_i,
1\le i\le 22$. Clearly, for each $P_i$ the triple of signs ordered
in accordance with the appearance of the products  in
(\ref{basicsigns}) is equal to $Sign(P_i)$ determined by the real
structure induced on $Y(\mathcal L)$ by the standard complex
conjugation in $\mathbb P^9_{w}$, $z_k\mapsto\bar z_k$ and
$u_\alpha\mapsto\bar u_\alpha$. (Any real structure on $Y(\mathcal
L)$ lifts to a real structure on $X(\mathcal L)$ and such a lift
is unique, cf. subsection \ref{Kl}.)

By Proposition \ref{klcam}, there are eight and only eight
distinct real structures $c_X$ which are lifts of $c_P$. Let show
that they all can be induced by a suitable diagonal real structure
on $\mathbb P^9_{w}$, where by a diagonal real structure on
$\mathbb P^9_{w}$ we mean a real structure given by
$z_k\mapsto\bar z_k$ and $u_\alpha\mapsto\epsilon_\alpha\bar
u_\alpha$ with $\epsilon_\alpha=\pm 1$. Note that such a real
structure $c_\epsilon$ preserves $Y(\mathcal L)$ if, and only if,
the equations (\ref{dependness}) are respected. In particular,
there are eight and only eight real diagonal structures which
preserve $X(\mathcal L)$ and they are determined by an arbitrary
choice of $\epsilon_\alpha$ with $\alpha=(1,0,0), (0,1,0)$, and
$(0,0,1)$. We denote by
$$
c_{
\epsilon_{(1,0,0)},\epsilon_{(0,1,0)},\epsilon_{(0,0,1)}}:X(\mathcal
L)\to X(\mathcal L)
$$
the real structures thus
obtained. Each of them
is a lift of $c_P$, since they all transform $z_k$ in $\bar z_k$.

As is easy to check, the sign-triple $Sign'_i=({sign}'_{i,1},
{sign}'_{i,2}, {sign}'_{i,3})$ of $P_i$ defined by
$c_{\epsilon_{(1,0,0)},\epsilon_{(0,1,0)},\epsilon_{(0,0,1)}}$ is
equal to $(\epsilon_1{sign}_{i,1}, \epsilon_2{sign}_{i,2}, \epsilon_3{sign}_{i,3}),$
which in its turn is equal to the triple of sings of the
homogeneous forms
\begin{equation} \label{signs}
\begin{array}{l}
\epsilon_{(1,0,0)}l_{(1,0,0)} l_{(1,1,0)} l_{(1,0,1)} l_{(1,1,1)}
,
\\
\epsilon_{(0,1,0)} l_{(0,1,0)} l_{(1,1,0)} l_{(0,1,1)} l_{(1,1,1)}
,
\\
\epsilon_{(0,0,1)}l_{(0,0,1)} l_{(1,0,1)} l_{(0,1,1)} l_{(1,1,1)}
.
\end{array}
\end{equation}
In what follows, a line
arrangement $\mathcal L$ equipped with one of these {\bf eight
sign-labelings} is called {\it equipped {\rm(}by signs{\rm)}}.

The sign-equipment of a (labelled) pure real Campedelli
arrangement contains a complete information on the real structure,
as the following proposition shows.

\begin{prop}\label{purerealclas}
Let Campedelli line arrangements $\mathcal L$
and $\mathcal L'$ be pure real with respect to
real structures $c_P:\mathbb P^2\to\mathbb P^2$ and
$c'_P:\mathbb P^2\to\mathbb P^2$.
A real structure $c: X(\mathcal L)\to X(\mathcal L)$ lifting
$c_P$ and a real structure
$c': X(\mathcal L')\to X(\mathcal L')$ lifting $c'_P$ are equivalent
if, and only if, there exist a homomorphism $h:(\mathbb Z/ 2\mathbb Z)^3
\to (\mathbb Z/ 2\mathbb Z)^3$ and a projective transformation
$H:\mathbb P^2\to\mathbb P^2$ such that
$$
c'_P\circ H=H\circ c_P, \quad \phi'\circ H_*=h\circ\phi
$$
(here $\phi: H_1(\mathbb P^2\setminus\mathcal L; \mathbb Z) \to
(\mathbb Z/ 2\mathbb Z)^3$ and $\phi': H_1(\mathbb
P^2\setminus\mathcal L'; \mathbb Z) \to (\mathbb Z/ 2\mathbb Z)^3$
are the labelings participating in definition of $\mathcal L$ and
$\mathcal L'$), and
$$
Sign'(H(P_i))=(-1)^{h(g_i)},
$$
where
$$
(-1)^{g_i}=Sign P_i.
$$
\end{prop}

\proof It follows from Theorem \ref{autcam} and the definition of
the sign-triples (recall that one sign-triple determines all the other).
\qed

\begin{prop}\label{eight}
The eight real structures $ c_{
\epsilon_{(1,0,0)},\epsilon_{(0,1,0)},\epsilon_{(0,0,1)}}$ are
distinct. If $\mathcal L$ has no triple points, these eight real
structures are the only reals structures of $X(\mathcal L)$.
\end{prop}

\proof There exist points in $Y(\mathcal L)$ where all the three
coordinates $z_0,z_1,$ $z_2$ are real and all the three
coordinates $u_{(1,0,0)},u_{(0,1,0)},u_{(0,0,1)}$ are nonzero. The
real structures
$c_{\epsilon_{(1,0,0)},\epsilon_{(0,1,0)},\epsilon_{(0,0,1)}}$
with different $(\epsilon_{(1,0,0)},$ $\epsilon_{(0,1,0)},$ $
\epsilon_{(0,0,1)})$ act differently on such a point. It implies
the first statement.

Now, assume that $\mathcal L$ has no triple points and consider
two real structures,
$c: X(\mathcal L)\to X(\mathcal L)$ lifting
$c_P$ and
$c': X(\mathcal L)\to X(\mathcal L)$ lifting $c'_P$.
Assume that $\mathcal L$ is pure real with respect to $c_P$.

Let show, first, that $\mathcal L$ is pure real with respect to $c'_P$
as well. Suppose that
$\mathcal L$ is mixed real with respect to $c'_P$. Then, $c'_P$
preserves three lines in $\mathcal L$ and two points which are intersections
of conjugated lines in $\mathcal L$. On the other hand, $c_P$ preserves
all the lines in $\mathcal L$. Therefore, $c_P\circ c'_P$ acts trivially
on three generic lines and two points outside them. Hence, $c_P\circ c'_P=
\text{id}$, so that $c_P= c'_P$. Which is a contradiction.

If $\mathcal L$ is pure real with respect to both $c'_P$ and
$c_P$, then the same argument implies that $c_P= c'_P$. Thus, $c$
and $c'$ differ by a Galois transformation. \qed

\begin{rem} As follows from
Propositions \ref{purerealclas} and \ref{eight}, if $\mathcal L$
is a pure real Campedelli arrangement without nontrivial projective
automorphisms, then the eight real structures
$c_{\epsilon_{(1,0,0)},\epsilon_{(0,1,0)},\epsilon_{(0,0,1)}}$
are non-equivalent to each other
and they represent all the real structures on $X(\mathcal L)$.
\end{rem}

\begin{lem} \label{positive}
For any choice of $\epsilon_{\alpha}$ with $\alpha=(1,0,0)$, $(0,1,0)$, and $(0,0,1)$,
the real point set $X_{\mathbb R}=\text{Fix}\,\, c$,
$c=c_{\epsilon_{(1,0,0)},\epsilon_{(0,1,0)},\epsilon_{(0,0,1)}}$, is
$$
X_{\mathbb R}= \bigcup_{Sign_i=(+,+,+)}f^{-1}(P_i),
$$
where $Sign_i$ are the sign-triples defined by $c$. \qed\
\end{lem}

Assume that $\mathcal L$ has no triple points (in fact,
one can treat in a similar way the
degenerate cases, but we do not need it). Let $P_{i_0}$ be a
$n$-gon. For each its side and for each its vertex, there is one and
only one polygon $P_i, i\ne i_0,$ intersecting $P_{i_0}$ along this
side or, respectively, at the vertex. Inspecting the sides and the
vertices along the border of $P_{i_0}$, we obtain a sequence of
polygons
$$(P_{i_1},P_{i_2}',\dots,P_{i_{2n-1}},P_{i_{2n}}'),$$ where $P'_- $
are the polygons adjacent to the vertices, and associate with it an
integer sequence
$A_{i_0}=(n_{i_1},n_{i_2}',\dots,n_{i_{2n-1}},n_{i_{2n}}')$, where
$n_{i_j}$ and $n_{i_j}'$ state for the number of sides of $P_{i_j}$
and, respectively, $P_{i_j}'$. The sequence $A_{i_0}$ is called the
{\it adjacency type} of $P_{i_0}$. The adjacency type is defined up
to cyclic permutation and reversing the order.

Let finally $\mathcal L$ be equipped by signs and let
$P_{i_1},\dots,P_{i_k}$ be the set of positive polygons. The
unordered collection $A(\mathcal L)=(A_{i_1},\dots, A_{i_k})$,
where $A_{i_j}$ is the adjacency type of $P_{i_j}$, is called the
{\it adjacency type of positive polygons}.

\begin{lem} \label{seven}
If $\mathcal L$ is a purely real Campedelli line arrangement
without triple points, then any its sign-equipment contains at
least seven different labels $Sign_i$.
\end{lem}

\proof The arrangement $\mathcal L$, as any arrangement without
triple points consisting of $\ge 5$ lines, defines at least five
triangles $P_i$. Through a simple counting of edges and cells, it
implies that in the case of seven lines there is a $n$-gon $P_i$
with $n\geq 5$.

If $P_i$ is a $\ge 6$-gon, then $P_i$ and the seven polygons
having a common side with $P_i$ have all different signs, as it
follows from the transition rule (\ref{transition}).

Let $P_i$ be a $5$-gon bounded by $L_{i_1}, L_{i_3}, L_{i_5},
L_{i_7},$ and $L_{i_9}$, and let
$(P_{i_1},P_{i_2}',\dots,P_{i_{9}},P_{i_{10}}')$ be its sequence
of adjacent polygons. As in the proof of Lemma \ref{gluing}, we
can assume (maybe, after renumbering of lines and a cyclic
permutation of adjacent polygons; note that a renumbering may
change the sign-equipment but preserve distinct the distinct
sign-triples) that $\alpha_{i_1}=(1,0,0)$, $\alpha_{i_3}=(0,1,0)$,
$\alpha_{i_5}=(0,0,1)$ and
$$\alpha_{i_7}, \alpha_{i_9}\in \{
(1,1,0),(1,0,1),(0,1,1),(1,1,1)\} .$$ By the transition rule
(\ref{transition}), the sign-triples of $P_i$ and its adjacent
polygons form the set $\{(-1)^a Sign_i$, $a\in A\}$, where $A=\{0,
\alpha_{i_1},$ $ \alpha_{i_3},$ $\alpha_{i_5},$ $\alpha_{i_7},$ $
\alpha_{i_9},$ $ \alpha_{i_1}+\alpha_{i_3},$
$\alpha_3+\alpha_{i_5},$ $ \alpha_{i_5}+\alpha_{i_7},$ $
\alpha_{i_7}+\alpha_{i_9},$ $\alpha_{i_9}+\alpha_{i_1}\}$. We have
$(0,0,0), (1,0,0),(0,1,0),(0,0,1),(1,1,0),(0,1,1)\in A$, that is,
$A$ consists of at least six elements. If $\alpha_7$ or $\alpha_9$
is equal to $(1,0,1)$ or $(1,1,1)$, then $A$ consists of at least
seven elements. Otherwise, if $\{ \alpha_7, \alpha_9\}= \{
(1,1,0),(0,1,1)\} $, again $A$ consists of at least seven
elements, since in this case $(1,1,0)+(0,1,1)=(1,0,1)\in A$. \qed

\begin{prop}\label{sevenreal}
Let $\mathcal L$ be a purely real Campedelli line arrangement
without triple points. For each real structure
$c_{\epsilon_{(1,0,0)},\epsilon_{(0,1,0)},\epsilon_{(0,0,1)}}$ on
$X=X(\mathcal L)$ except, possibly, one, its real points set is
non-empty.
\end{prop}

\proof It follows from Lemma \ref{positive} and Lemma \ref{seven}.
\qed

\subsection{Mixed real Campedelli line arrangements} \label{mixed}
\label{c-c} Let $\mathcal L$ be a Cam\-pe\-del\-li line arrangement
$\mathcal L= \cup L_{\alpha}$ which is mixed real with respect to a
real structure $c_P: \mathbb P^2\to\mathbb P^2$. Choose homogeneous
coordinates $(z_0,z_1,z_2)$ in $\mathbb P^2$ such that $c_P$ turns
in
$$
c_P(z_0,z_1,z_2)=(\overline z_0,\overline z_1,\overline z_2).$$
Then, up to a renumbering and a real projective transformation, the
lines $L_{(1,1,0)}$, $L_{(1,1,1)}$, and $L_{(0,0,1)}$ are given by
equations $z_0=0$, $z_1=0$, and $z_2=0$, while the lines
$L_{(1,0,0)}$, $L_{(0,1,0)}$, $L_{(1,0,1)}$, and $L_{(0,1,1)}$ are
given by equations
$$a_{\alpha,0}z_0+a_{\alpha,1}z_1+a_{\alpha,2}z_2=0,$$
where $a_{(1,0,0),j}=\overline a_{(0,1,0),j}$ and
$a_{(1,0,1),j}=\overline a_{(0,1,1),j}$ for any $j=0,1,2$ (cf., the
proof of Proposition \ref{klcam}).

As above, consider the set $\mathbb R\mathbb P^2=\{ (z_0:z_1:z_2)
\mid \, z_i\in \mathbb R\}$ of real points of $\mathbb P^2$.  A
mixed real Campedelli line arrangement $\mathcal L$ intersect
$\mathbb R\mathbb P^2$ along three distinct real lines
$L_{\alpha,\mathbb R}=L_{\alpha}\cap \mathbb R\mathbb P^2$,
$\alpha=(0,0,1), (1,1,1)$, and $(1,1,0)$, and at two distinct real
points $p_1=L_{(1,0,0)}\cap L_{(0,1,0)}$ and $p_2=L_{(1,0,1)}\cap
L_{(0,1,1)}$. We call the points $p_1$ and $p_2$ the {\it
vertices} of $\mathcal L$. The vertices can not belong to
$L_{1,1,0}$, but it may happen that one of them (or both together)
belong to $L_{1,1,1}\cup L_{0,0,1}$. There is a renumbering which
exchange $p_1$ and $p_2$.

Denote by
$l_{\alpha}=a_{\alpha,0}z_0+a_{\alpha,1}z_1+a_{\alpha,2}z_2$,
$\alpha\in G\setminus\{0\}$, the above linear forms defining
$L_{\alpha}$. Put $q_1=l_{(1,0,0)}l_{(0,1,0)}$ and
$q_2=l_{(1,0,1)}l_{(01,1)}$. Note that $q_1$ and $q_2$ have real
coefficients. Moreover, $q_1\ge 0$ and $q_2\ge 0$ at each point of
$\mathbb RP^2$.

The Campedelli surface $X=X(\mathcal L)$ is given in $\mathbb P^9_w$
by equations

\begin{equation} \label{equations4}
\begin{array}{l}
u_{(1,0,0)}^2=l_{(1,0,0)}l_{(1,0,1)}z_0z_1, \\
u_{(0,1,0)}^2=l_{(0,1,0)}l_{(0,1,1)}z_0z_1, \\
u_{(0,0,1)}^2=q_2z_1z_2, \\
u_{(1,1,0)}^2=q_1q_2, \\
u_{(1,0,1)}^2=l_{(1,0,0)}l_{(0,1,1)}z_0z_2, \\
u_{(0,1,1)}^2=l_{(0,1,0)} l_{(1,0,1)}z_0z_2, \\
u_{(1,1,1)}^2=q_1z_1z_2.
\end{array}
\end{equation}
It inherits a real structure $c_{++}:X\to C$ from the real structure
on $\mathbb P^9_w$ defined by $z_k\mapsto \bar z_k$ and
$u_{(i,j,k)}\mapsto \bar u_{(j,i,k)}$.

In accordance with Proposition \ref{klcam}, there are three more
real structures on $X$ (only three, if the arrangement has no a
nontrivial projective automorphism) which are obtained from
$c_{++}$ by composing it with Galois actions. Namely, they are
\begin{equation} \label{structure4}
c_{-+}=g_{(1,0,0)}c_{++}g_{(1,0,0)}, c_{+-}=g_{(0,0,1)}c_{++},
\quad\text{and}\quad c_{--}=g_{(0,0,1)}c_{-+},
\end{equation}
where $g_{(1,0,0)}$, $g_{(0,0,1)}\in \text{Gal}(X/\mathbb P^2)$ are
defined as follows:
$$
g_{(1,0,0)}u_{(i,j,k)}=(-1)^iu_{(i,j,k)}, \,\,
g_{(0,0,1)}u_{(i,j,k)}=(-1)^ku_{(i,j,k)}.
$$
In particular, one can notice that up to conjugation by
automorphisms of $X$ this list of four real structures reduces to
two conjugacy classes represented, respectively, by $c_+=c_{++}$ and
$c_-=c_{+-}$.

Let assume that $\mathcal L$ has no triple points and subdivide
such arrangements in three types. The lines $L_{\alpha, \mathbb
R}$, $\alpha=(0,0,1), (1,1,1), (1,1,0)$, divide $\mathbb R\mathbb
P^2$ into four triangles $P_i$, $i=1,\dots,4$, as it is depicted
in Fig.2, where the axe $x=0$ is the line $L_{(1,1,1),\mathbb R}$,
the axe $y=0$ is the line $L_{(0,0,1),\mathbb R}$, while the line
$L_{(1,1,0),\mathbb R}$ is put at infinity.
Using renumberings which transform $(1,0,0)$ in $(1,0,0)$,
$(0,1,0)$ in $(0,1,0)$, and $(0,0,1)$ in $(1,1,1)$ together with
linear transformations $x\mapsto \pm x$, $y\mapsto \pm y$ , we can
and will assume that $p_1\in P_1$ and
$p_2$ belongs either to $P_1$ ({\bf Type I}), or to $P_2$ ({\bf Type
II}), or to $P_3$ (\bf Type III}).

\begin{picture}(300,130)
\put(90,65){\vector(1,0){120}}\put(150,25){\vector(0,1){100}}
\put(170,80){$P_1$}\put(170,35){$P_2$}\put(120,80){$P_4$}\put(120,35){$P_3$}
\put(154,118){$y$}\put(207,57){$x$} \put(140,5){$\text{Fig. 2}$}
\end{picture}

Such a normalization makes the products $l_{(1,1,0)} l_{(1,1,1)}=
z_0z_1$ and $l_{(1,1,0)}l_{(0,0,1)}=z_0z_2$ to be positive on
$P_1$ (and on $P_3)$ and, in particular, fixes a choice of $c_+$.
Under this convention,  $c_-$ becomes the real structure induced
by $z_k\mapsto \bar z_k$ and $u_{(i,j,k)}\mapsto \bar u_{(j,i,k)}$
on the copy of $X$ which is given by
\begin{equation} \label{equations5}
\begin{array}{l}
u_{(1,0,0)}^2=l_{(1,0,0)})l_{(1,0,1)}v_1,
\\
u_{(0,1,0)}^2=l_{(0,1,0)}l_{(0,1,1)}v_1,
\\
u_{(0,0,1)}^2= -q_{2}v_1v_2.
\end{array}
\end{equation}

\begin{lem} \label{ll} Let
$\mathcal L$ be a mixed real Campedelli line arrangement without
triple points. Suppose that $P_i$ and $c_\pm$ are labelled as
above. Then, for any $i=1,\dots,4$,
\begin{itemize} \item[($i$)] $f^{-1}(P_i)$
is a disjoint union $P_{i,1}\cup P_{i,2}$ of two connected
non-orientable two-manifolds, \item[($ii$)] the Euler
characteristic of $P_{i,j}, j=1,2,$ is equal to $1-2n$, where $n$
is the number of vertices $\{ p_1,p_2\}$ belonging to $P_i$,
\item[($iii$)] the real point set $X_{\mathbb R}=\text{Fix}\, c$,
$c=c_\pm$, is
$$
X_{\mathbb R}= P_{i,1}\cup  P_{i+2,1},
$$
where $i=1$ if $c=c_+$ and
$i=2$ if $c=c_-$.
\end{itemize}
\end{lem}
\proof It is similar to the proof of Lemma \ref{lll}. The only
difference is that here inside $P_i$ we have vertices $p_1,p_2$
which are (simple) branching points of the projection $\widetilde
P_{i,j}\to P_i$.\qed

\begin{rem}
A Campedelli line arrangement can be purely real with respect to one real structure
and mixed real with respect to another one. More precisely,
a Campedelli line arrangement $\mathcal L$ is simultaneously purely
real and mixed real if and only if (maybe after
renumbering of the lines) there are coordinates $(z_0,z_1,z_2)$ in
$\mathbb P^2$ such that: $z_i=0$, $i=0,1,2$, is an equation
respectively of $L_{(1,1,0)}$, $L_{(1,1,1)}$, $L_{(0,0,1)}$; the
lines $L_{(1,0,0)}$ and $L_{(0,1,0)}$ are given by $a_1z_1
+(a_0z_0\pm z_2)=0$, and the lines $L_{(1,0,1)}$ and $L_{(0,1,1)}$
are given by $b_1z_1+( b_0z_0\pm z_2)=0$ for some non-zero $a_0$,
$a_1$, $b_0$, $b_1\in \mathbb R$.
\end{rem}

\section{Diffeomorphisms and deformations of real Campedelli surfaces}
\label{difanddef}

\subsection{Deformation versus smoothing of $A_1$-points}
By a real Morse-Lefschetz perturbation of a real surface with
$A_1$-singu\-larities we mean a complex three-manifold $Z$ with a
real structure $c:Z\to Z$ equipped with a proper holomorphic map
$f$ from $Z$ to the unit disc $D\subset{\mathbb C}$ respecting the
real structures on $Z$ and $D\subset \mathbb C$ and such that: all
the fibers of $f$, except the fiber over $0$, are (compact)
nonsingular surfaces; the fiber over $0$ contains only isolated
singular points $O_1, \dots, O_k$, and the quadratic form of $f$
at each of the singular points is non-degenerate. The fibers
$f^{-1}(t)$ are denoted by $X_t$, so that the singular fiber
$f^{-1}(0)$ is denoted by $X_0$. The real structure $c:X_0\to X_0$
lifts to a unique real structure $c:\tilde X_0\to \tilde X_0$
where $\tilde X_0$ is the minimal desingularization of $X_0$.
According the definition of the deformation equivalence of real
surfaces, for all $t\in{\mathbb R}, t\ne 0,$ of the same sign the
real surfaces $(X_t,c)$ are of the same real deformation type. If
$O_j, 1\le j\le k,$ is real then we pick a small (Milnor) ball
$B_j\subset Z$ around $O_j$ and, for every small real $t\ne 0$,
speak on the local Euler characteristic of $X_{t,\mathbb R}$ which
means the Euler characteristic of the intersection of the real
part of $X_t$ with $B_j$.

Such Morse-Lefschetz perturbations are provided by triangle moves
of real Campedelli line arrangements, see subsection
\ref{reversing}.

\begin{lem}\label{morse} Let $(Z,f,c)$ be a real Morse-Lefschetz
perturbation of a real surface with $A_1$-singularities. If for
$t'\ne 0$ of certain sign, at each singular point $O_j\in X_0$
which is real the local Euler characteristic of
$X_{t,\mathbb R}$ is $0$, then $(X_{t'},c)$ is real deformation
equivalent to $(\tilde X_0, c)$.
\end{lem}

\proof Introduce an auxiliary real one-parametric family by making
the base change which substitutes $u^2$ instead of $t$ if $t'$ is
positive, and $-u^2$ otherwise. The total space of this family has
$A_1$-singularities at $O_1,\dots, O_k\in X_0$ and it has no any
other singular point. Blowing up the total space at the
$A_1$-singularities we respect the real structure and replace each
of the singular points by a quadric and resolve both the singular
points of the family and the singular points of $X_0$. At each
point $O_j$ which is real the blown-up quadric is real, and the
two families of generating lines on this real quadric are real if,
and only if, the local Euler characteristic of $X_{t,\mathbb R}$
with $t=t'$ is $0$. Pick a real family of lines at
each of real $O_j$ and conjugated families of lines at each pair
of conjugated $O_j$. As is known, a contraction of any family of
lines gives a smooth family. The contraction of the chosen
families is real and thus provides a real deformation equivalence
between $(\tilde X_0, c)$ and $(X_{t'},c)$. \qed

\begin{rem}\label{ccA1} If $(Z,f,c)$
is a Morse-Lefschetz perturbation of a real surface with complex
conjugated {\rm (}non real{\rm )} $A_1$-singularities, then all
$(X_{t},c)$ with real $t\ne 0$ are real deformation equivalent to
each other.
\end{rem}

\subsection
{Triangle moves}\label{reversing}
Let $\mathcal L$
be an equipped purely real Cam\-pe\-del\-li line arrangement, see
subsection \ref{purely}, and let $P_{i_0}\subset \mathbb R\mathbb
P^2$ be a triangle of $\mathcal L$ whose sides are $L_{\alpha_1}$,
$L_{\alpha_2}$, $L_{\alpha_3}$. A modification depicted in Fig. 3
which turns $\mathcal L$ into an equipped purely  real Campedelli
line arrangement $\mathcal{L}'$ is called the {\it reversing of
triangle $P_{i_0}$} or {\it a triangle move}. By definition, the
sign-triples $\text{Sign}_{i}'=\text{Sign}(P_{i}')$ with $i\neq
i_0$ coincide with $\text{Sign}_{i}=\text{Sign}(P_{i})$, while, in
accordance with the transition rule (\ref{transition}),
$\text{Sign}_{i_0}'=(\text{sign}_{i_0,1}',\text{sign}_{i_0,2}',\text{sign}_{i_0,3}')$
is determined by
$\text{Sign}_{i_0}=(\text{sign}_{i_0,1},\text{sign}_{i_0,2},\text{sign}_{i_0,3})$
as follows:
$$\text{sign}_{i_0,j}'=(-1)^{a_j}\text{sign}_{i_0,j},$$
where $(a_1,a_2,a_3)=\alpha_1+\alpha_2+\alpha_3$.

\begin{picture}(380,190)
\qbezier(0,110)(50,110)(100,110) \qbezier(40,60)(40,120)(40,180)
\qbezier(10,170)(50,130)(100,80) \put(45,30){$\mathcal L$}
\put(105,110){$\leftrightsquigarrow$} \put(45,115){$P_{i_0}$}
\put(41,172){$L_{\alpha_1}$}\put(1,114){$L_{\alpha_2}$}
\put(90,74){$L_{\alpha_3}$}

\put(25,165){$P_{i_2}$}
\put(20,130){$P_{i_3}$}\put(20,85){$P_{i_4}$}
\put(60,85){$P_{i_5}$}\put(90,95){$P_{i_6}$}
\put(65,135){$P_{i_1}$}

\qbezier(130,120)(180,120)(230,120)
\qbezier(180,60)(180,120)(180,180)
\qbezier(140,160)(180,120)(220,80)
\put(175,30){$\mathcal{L}^0$}\put(235,110){$\leftrightsquigarrow$}
\put(181,172){$L_{\alpha_1}$}\put(210,74){$L_{\alpha_3}$}
\put(133,109){$L_{\alpha_2}$}\put(155,160){$P_{i_2}$}
\put(150,130){$P_{i_3}$}\put(150,85){$P_{i_4}$}
\put(190,85){$P_{i_5}$}\put(210,105){$P_{i_6}$}
\put(195,140){$P_{i_1}$}

\put(336,119){$L_{\alpha_2}$}\put(237,164){$L_{\alpha_3}$}
\qbezier(250,130)(290,130)(350,130)
\qbezier(310,60)(310,120)(310,180)
\qbezier(250,160)(300,110)(340,70) \put(310,30){$\mathcal{L}'$}
\put(296,117){$P'_{i_0}$} \put(311,172){$L_{\alpha_1}$}
\put(275,155){$P'_{i_2}$} \put(255,135){$P'_{i_3}$}
\put(275,90){$P'_{i_4}$}\put(315,75){$P'_{i_5}$}
\put(330,95){$P'_{i_6}$} \put(325,150){$P'_{i_1}$}

\put(150,5){$\text{Fig. 3}$}
\end{picture}

\begin{rem}\label{rev} {\rm If the sides of $P_{i_0}$ are linear dependent,
then: $\text{Sign}_{i_j}=\text{Sign}_{i_{j+3}}$ for $j=1,2,3$;
$\text{Sign}_{i_0}'=\text{Sign}_{i_0}$; and $\text{Sign}_{i_1}$,
$\text{Sign}_{i_2}$, and $\text{Sign}_{i_3}$ are pairwise
distinct. In the case of linear dependent sides, $\mathcal{L}^0$
is not a Campedelli arrangement.}
\end{rem}

If the sides of $P_{i_0}$ are linear independent, all the triples
$\text{Sign}_{i_j}$, $j=0,1,\dots,6$, are pairwise distinct and
$\text{Sign}_{i_0}'$ is the complementary element in the set of
all triples of signs. In the case of linear independent sides,
$\mathcal{L}^0$ is a Campedelli line arrangement and the canonical
model $X(\mathcal{L}^0)$ of a Campedelli surface $X^0$ has two
$A_1$-singular points over the triple point (the point to which
$P_{i_0}$ degenerates). The local Euler characteristic of
$X(\mathcal L)$ at these singular points is $0$ if, and only if,
$$
(+,+,+)\notin\{Sing_{i_k}\}_{k=0,2,4,6}.
$$
Respectively, the local Euler characteristic of $X(\mathcal L')$ at these singular points is
$0$ if, and only if, $(+,+,+)\notin\{Sing'_{i_k}\}_{k=0,2,4,6}$.
The last condition is equivalent to
$(+,+,+)\in\{Sing_{i_k}\}_{k=0,2,4,6}$;
in particular, if the local Euler characteristic is equal to $0$
for one of $\mathcal L$ and $\mathcal L'$, it is not equal to $0$
for the other one, and vise versa.

\subsection{Reduction to generic deformations}

\begin{lem}\label{tomorse} Suppose that $(Z,f,c)$ is a real deformation
such that all the fibers except $X_0$ have nonsingular canonical
models, while the canonical model of $X_0$ is a surface with
$A_1$-singularities. Then, at each singular point $O_j\in Z$ which
is real the local Euler characteristic of
$X_{t,\mathbb R}$, $t\ne 0$ is $0$.
\end{lem}

\proof The deformation $(Z,f,c)$ is a simultaneous resolution of
the singularities of the family constituted of the canonical
models $X^{can}_t$ of $X_t$ and regarded over the same base.
Hence, for each small real $t$ the local Euler characteristics of
$X_{t,\mathbb R}$ coincide with the local Euler characteristics of
the resolutions of the singular points. The latter characteristics
are $0$ in the case of $A_1$-singularities, whatever are the real
forms of the singularities.\qed

\begin{lem}\label{realcanon} Let $(Z,f,c)$ be a real deformation
of Campedelli surfaces. For any real $t'\in D$, there exist a real
neighborhood $U\subset D$ of $t'$ and a real family $\mathcal L_t,
t\in U,$ of Campedelli line arrangements in a real projective
plane $(\mathbb P^2; c_P)$ such that $X_t=X(\mathcal L_t)$ and
$c_t=c\vert_{X_t}$ are lifts of $c_P$.
\end{lem}

\proof Consider the relative bi-canonical bundle $2K\vert_{Z/D}$.
Its restriction to any fiber $X_t$ is the bi-canonical bundle of
$X_t$. The space of sections of such a restriction is of dimension three,
and the sections determine a finite map to $\mathbb P^2$ representing
$X_t$ as $X(\mathcal L_t)$, where $\mathcal L_t$ is the branching locus of
this map, see the proof of
Theorem \ref{autcam}. Since the space of sections is of constant
dimension, any three sections generating
the bi-canonical bundle of $X_{t'}$ extend to three sections
generating  $2K\vert_{Z/D}$ at least over a small neighborhood of
$t'$. By theorem \ref{realthm}, the three sections of the bi-canonical
bundle of $X_{t'}$ can be chosen real with respect to a real
structure $c_P$ of $\mathbb P^2$, and then it remains to average
their extensions by $c$ and pick a sufficiently small equivariant
neighborhood of $t'$. \qed

\begin{prop}\label{defeq} Let $(X_1,c_1)$ and $(X_2,c_2)$ be two
deformation equivalent real Campedelli surfaces associated
respectively with Cam\-pe\-del\-li line arrangements $\mathcal
L_1$ and $\mathcal L_2$. If $\mathcal L_1$ is purely real, then
$\mathcal L_2$ is also purely real. If they are purely real and
have no triple points, then their sign-equipments in $\mathbb
RP^2$ are homeomorphic, so that, in particular, $\mathcal L_1$ and
$\mathcal L_2$ have the same type and the same adjacency type of
positive polygons.
\end{prop}

Note that this statement implies that  a deformation can not
provide any triangle move.

\proof By Lemma \ref{realcanon}, a chain of real deformations
connecting $(X_1,c_1)$ and $(X_2,c_2)$ results in a chain of real
families of Campedelli line arrangements $\mathcal L_t$. We look
at $\mathcal L_t$ with real values of $t$. It gives a chain of
real Campedelli line arrangements connecting $\mathcal L_1$ and
$\mathcal L_2$. Campedelli line arrangements have at worse triple
points. Therefore, the number of real lines in an arrangement is
not changing in a chain of real deformations. It proves the first
statement.

Now assume that $\mathcal L_1$ and $\mathcal L_2$ are purely real
and have no triple points. The triple points on intermediate
arrangements $\mathcal L_t$ appear and disappear independently.
Their appearance and disappearance befalls by triangle half-moves:
contracting and reappearing of triangles like in Fig. 3. The
half-move provided by a reappearing triangle should turn back the
local combinatorial structure and the local sign-equipment, since
according to Lemma \ref{tomorse} the local input to the Euler
characteristic of the real part should be $0$ for both types of
half-moves, while, as we observed already in subsection
\ref{reversing}, such an input due to a contracting triangle
$P_{i_0}$ (or, respectively, to a turning back triangle
$P_{i_0}'$) is $0$ if, and only if,
$(+,+,+)\notin\{Sing_{i_k}\}_{k=0,2,4,6}$ (respectively,
$(+,+,+)\notin\{Sing'_{i_k}\}_{k=0,2,4,6}$). Finally, it implies
that replacing the straight lines $L_{i,\mathbb R}$ by suitable
flexible lines one can connect $\mathcal L_1$ and $\mathcal L_2$
by a continuous family of equipped flexible configurations in
$\mathbb R\mathbb P^2$ without triple points.\qed

\subsection{Smoothing of $T(-4)$ singularities.}\label{Tsing}
By {\it a real smoothing } of a real surface $(M,c)$ we mean any
real fiber of a real flat family of surfaces $Z\to D$ over the
unit disc $D$ (where the real structure on $D$ is given by the
usual complex conjugation) such that $(X_0,c)=(M,c)$ and $X_t$ is
nonsingular for any $t\in D, t\ne 0$. A singular point of a
surface is called {\it $T(-4)$-singularity} if its germ is
isomorphic to the $(\mathbb Z/2\mathbb Z)^2$-Galois covering of
the germ $(\mathbb C^2, 0)$ branched in three lines
$L_{\alpha_1}\cup L_{\alpha_2}\cup L_{\alpha_3}$ through $0$ with
a $(\mathbb Z/2\mathbb Z)^2$ labeling $\{\alpha_i\}_{i=1,2,3}$
such that $\alpha_1+\alpha_2+\alpha_3=0$. We speak on {\it a real
surface with non real $T(-4)$-singularities}, if all the singular
points of the surface are $T(-4)$-singularities and neither of the
singular points is real.

\begin{thm} \label{throughtriple}
Any two real smoothings $(M_1,c)$ and $(M_2,c)$ of a real surface
$(M,c)$ with non real $\,T(-4)$ singularities
 have diffeomorphic real structures.
\end{thm}

\proof The pairs $(M_1,c)$ and $(M_2,c)$ are obtained from $(M,c)$
by removing $c$-invariant Milnor neighborhoods $U_j\cup c(U_j)$ of
the each pair of conjugated singularities followed by a
$c$-invariant gluing of some standard pieces $N_j\cup N_{\bar j}$,
$N_j=(N,j)$ and $N_{\bar j}=(N,\bar j)$, instead of $U_j\cup
c(U_j)$ by means of some boundary diffeomorphisms $\phi_j:\partial
N\to \partial U_j, \phi_{\bar j}:\partial N\to \partial c(U_j)$
such that $c\circ \phi_j=\phi_{\bar j}$ (so that $c$ acts on
$N_j\cup N_{\bar j}$ by $(x,j)\mapsto (x,\bar j)$. As is shown,
for example, in \cite{Man}, the result of gluing of the half of
these pieces, say $\cup N_j$, gives diffeomorphic four-manifolds
$M_1\setminus\bigcup_j N_{\bar j}$ and $M_2\setminus\bigcup_j
N_{\bar j}$ (in fact, $\partial N$ is a lens space $L(4,1)$ and
the existence of such a diffeomorphism follows from a
corresponding Bonahon theorem, see \cite{Bon}). Now, it remains to
extend such a diffeomorphism $\Phi$ to $M_1\to M_2$ by symmetry,
that is by taking $\Phi(x)=(c\circ \Phi)(x)$ for each $x\in
(N,\bar j)$.\qed

The following corollary is straightforward.

\begin{cor}\label{rev-depend}
Let $\mathcal L$ and $\mathcal L'$ be two equipped real Campedelli
line arrangements related by a triangle move reversing a triangle
$P_{i_0}$. Suppose that the sides of $P_{i_0}$ are linear
dependent and that $\text{Sign}_{i_j}\neq (+,+,+)$ for
$j=0,1,\dots,6$, where $P_{i_j}$,
$j=1,\dots,6$, are the polygons adjacent to $P_{i_0}$. Then the
real Campedelli surfaces $X(\mathcal L)$ and $X(\mathcal{L}')$
have diffeomorphic real structures. \qed
\end{cor}

\subsection{
Classification of mixed real Cam\-pe\-del\-li line arrangements up
to real deformations}\label{defsubs} Let $\mathcal L$ be a
Campedelli line arrangement which is mixed real with respect to a
real structure $c_P:\mathbb P^2\to\mathbb P^2$. We say that a real
Campedelli surface $(X,c_X)$, where $X=X(\mathcal L)$ and $c_X$ is
a lift of $c_P$, has the {\it type} $J_{\pm}$, where $J=I$, $II$,
or $III$, if: $\mathcal L$ is without triple points; it has the
type $J$; and $c_X=c_{\pm}$ (see subsection \ref{mixed} for
notation $I,II,III$ and $c_\pm$).

\begin{thm} \label{def-c-c}
There are exactly five types of deformation non-equi\-valent real
Campedelli surfaces $(X,c)$ associated with mixed real
Campe\-delli line arrangements. They are represented by
arrangements of types $I_{\pm}$, $II_+$, and $III_{\pm}$.
\end{thm}

\proof  According to Proposition \ref{defeq}, $X(\mathcal L)$ and
$X(\mathcal L')$ are not deformation equivalent if $\mathcal L$ is
a purely real Campedelli arrangement and $\mathcal L'$ is a mixed
real one. By Lemma \ref{morse}, if $\mathcal L$ has triple points,
the surface $(X,c_X)$ is real deformation equivalent to a surface
associated with a mixed real Campedelli line arrangement without
triple points. Therefore, there exist at most six types of
deformation non-equivalent real Campedelli surfaces $(X,c)$,
associated with mixed real Campedelli line arrangements, namely:
$I_{\pm}$, $II_{\pm}$, and $III_{\pm}$.

To distinguish them, notice that a real deformation of
a real Campe\-delli surface $(X,c_X)$ is simultaneously a
$H=\text{Kl}(X,c_P)$-defor\-ma\-tion, in a sense that not only the action of
$c_X$ but the action of the whole group $H$ extends to the total
space of the deformation. Moreover, since the Galois group
$G\subset H$ preserves each fiber of the deformation, the real
deformation of $(X,c_X)$ is simultaneously a real deformation for
each of the other real structures contained in $H$.

In the case of
mixed real Campedelli arrangements, $H$ is a quaternion group
(see \ref{klcam}), it contains four distinct real structures, and
they split in two conjugacy classes $c_\pm$ (see subsection
\ref{mixed}). Finally, the topological type of the unordered pair
of two-manifolds $(\text{Fix} c_+, \text{Fix}c_-)$ is invariant
under real deformations of $(X,c_X)$. Lemma \ref{ll} implies
this invariant distinguishes the cases $I_{\pm}$, $II_+$,
and $III_{\pm}$.

To finish the proof, let us show that the types $II_+$ and $II_-$
are deformation equivalent. Up to deformation equivalence, we can
assume that the vertices $p_1$ and $p_2$ of an arrangement
$\mathcal L$ of type $II$ have projective coordinates $(1,1,1)$
and, respectively, $(1,1,-1)$; moreover, we can assume that
$l_{(1,0,0)}l_{(0,1,0)}=(z_1-z_0)^2+(z_2-z_0)^2$ and,
respectively, $l_{(1,0,1)}l_{(0,1,1)}=(z_1-z_0)^2+(z_2+z_0)^2$.
Then, the diagonal transformation $z_0\mapsto z_0$, $z_1\mapsto
z_1$, $z_2\mapsto - z_2$ gives rise (see equations
(\ref{equations4})) to an equivalence between the real structures
$c_-$ and $c_+$ after renumbering $(1,0,0)\mapsto (1,0,1)$,
$(0,1,0)\mapsto (0,1,1)$, and $(0,0,1)\mapsto (0,0,1)$. \qed

\section{"Dif $\neq$ Def"}
\subsection{Example with a pair of purely real
arrangements}\label{7-gone}
\begin{ex} \label{example}  Two real Campedelli surfaces which have diffeomorphic
real structures but which are not real deformation equivalent.
\end{ex}

Let $\mathcal L$ be a purely real Campedelli line arrangement
defined by the sides of a $7$-gone $P_1$, that is an arrangement
of type $(7,14,0,0,1)$. Label the sides in a way that three
consecutive sides of $P_1$ be labelled by $(1,0,0), (1,1,0)$, and
$(0,1,0)$. Then, sign the triangle $P_0$ having a common side with
$P_1$ along $L_{(1,1,0)}$ by $(-,-,-)$ and extend this choice to a
sign-equipment of $\mathcal L$ following the transition rule
(\ref{transition}), see a fragment in  Fig. 4.

\begin{picture}(380,190)
\qbezier(100,110)(150,110)(200,110)
\qbezier(140,60)(140,120)(140,180)
\qbezier(110,170)(150,130)(200,80) \put(145,30){$\mathcal L$}
\put(142,120){$P_0$}
\qbezier(143,114)(145,114)(146,114)\qbezier(148,114)(150,114)(151,114)
\qbezier(153,114)(155,114)(156,114)
\put(141,172){$L_{(1,0,0)}$}\put(73,114){$L_{(0,1,0)}$}
\put(190,74){$L_{(1,1,0)}$} \put(105,175){$-+-$}
\put(100,130){$+--$}\put(100,85){$++-$}
\put(150,85){$-+-$}\put(190,95){$+--$} \put(165,135){$++-$}

\put(135,5){$\text{Fig. 4}$}
\end{picture}

Let $\mathcal L'$ be a sign-equipped arrangement obtained by the
triangle move reversing the triangle $P_0$. This arrangement is of
type type $(7,13,1,1,0)$. By Corollary \ref{rev-depend}, the real
Campedelli surfaces $X(\mathcal L)$ and $X(\mathcal{L}')$ have
diffeomorphic real structures, and by Proposition \ref{defeq} they
are not deformation equivalent.

\subsection{Example with eight purely real
arrangements}\label{eight}

\begin{ex} \label{main}
Eight real Campedelli surfaces $(X_1,c_1)$, $\dots$, $(X_8,c_8)$
which have diffeomorphic to each other real structures and which
are pairwise non-deformation equivalent.
\end{ex}

As in Example \ref{7-gone}, we search for equipped purely real
Campedelli line arrangements $\mathcal{L}_i$, $i=1,\dots,$, such
that, first, they are related be sequences of triangle moves
reversing non-positive triangles with linear dependent sides, and,
second, they differ by their types or the adjacency types of their
positive polygons. Then, by Theorem \ref{throughtriple}, the real
Campedelli surfaces $X(\mathcal{L}_i)$ have diffeomorphic real
structures, and by Corollary \ref{defeq}, they are pairwise
non-deformation equivalent.

\begin{picture}(350,450)
\qbezier(0,350)(130,350)(370,350)
\qbezier(68,20)(213,220)(358,420)
\qbezier(35,220)(200,310)(365,400)
\qbezier(250,20)(205,230)(160,440)
\qbezier(108,20)(139,230)(170,440)
\qbezier(40,440)(150,230)(260,20)
\qbezier(20,440)(180,285)(340,130) \put(5,420){\small
$L_{(0,1,1)}$} \put(130,435){\small $L_{(1,0,1)}$}
\put(33,240){\small $L_{(1,1,0)}$} \put(12,353){\small
$L_{(1,0,0)}$} \put(175,433){\small $L_{(0,0,1)}$}
\put(53,46){\small $L_{(0,1,0)}$} \put(43,435){\small
$L_{(1,1,1)}$} \put(91,353){\small{$P_1$}}
\put(177,289){\small{$P_2$}} \put(163,353){\small{$P_3$}}
\put(296,353){\small{$P_4$}} \put(133,265){\small{$P_5$}}
\put(194,183){\small{$P_6$}}
\put(30,385){\small{$P_7$}}\put(253,123){\small{$P_8$}}
\put(123,385){\small{$P_9$}} \put(170,67){\small{$P_{10}$}}
\put(243,385){\small{$P_{11}$}} \put(88,170){\small{$P_{12}$}}
\put(50,315){\small{$P_{13}$}}\put(120,315){\small{$P_{14}$}}
\put(133,333){\small{$P_{15}$}} \put(162,325){\small{$P_{16}$}}
\put(195,333){\small{$P_{17}$}} \put(248,315){\small{$P_{18}$}}
\put(150,295){\small{$P_{19}$}}\put(173,243){\small{$P_{20}$}}
\put(208,238){\small{$P_{21}$}} \put(146,187){\small{$P_{22}$}}
\put(160,-5){$\text{Fig. 5}$}
\end{picture}

To construct such arrangements, start from a purely real
Campedelli line arrangement $\mathcal{L}_{(0,0,0,0)}$ of type
$(11,5,5,1,0)$ depicted in Fig 5. This arrangement has six
pairwise disjoint triangles $P_1,\dots, P_6$. Each of them has
linear dependent sides. The number of sides for each of
$P_7,\dots,P_{22}$ is given in Table $(0,0,0,0)_1$.

Sign the triangle $P_1$  by $(+,+,+)$ and extend this choice to a
sign-equipment of $\mathcal L$ following the transition rule
(\ref{transition}). Then, as is easy to check,
$\mathcal{L}_{(0,0,0,0)}$ has only two positive polygons, namely
$P_1$ and $P_2$. \vspace{0.1cm}
\begin{center}
\begin{tabular}{|c|c|c|c|c|c|c|c| 
} \hline
   $P_7$ & $P_8$ & $P_9$ & $P_{10}$ & $P_{11}$ & $P_{12}$ & $P_{13}$ & $P_{14}$ \\
   \hline
 4 & 4 & 5 & 4 & 5 & 4 & 4
& 5 \\ \hline \hline
 $P_{15}$
   & $P_{16}$ & $P_{17}$ & $P_{18}$ & $P_{19}$ & $P_{20}$ & $P_{21}$ &
   $P_{22}$
 \\ \hline 3 & 5 & 3 & 5 & 3 & 6 & 3 & 3 \\
   \hline
\end{tabular}
\end{center} \vspace{0.1cm}

\begin{center}
Table $(0,0,0,0)_1$
\end{center}
\vspace{0.1cm}

To insure a possibility to perform independent triangle moves
reversing the four triangles $P_3,\dots, P_6$ it is sufficient to
consider $\mathcal{L}_{(0,0,0,0)}$ as a perturbation of a
degenerate configuration shown in Fig. 6. Now, it remains to
select the moves and to count for each configuration its type and
the adjacency type of its positive polygons.

\begin{picture}(350,250)
\qbezier(0,200)(130,200)(350,200)
\qbezier(40,122)(100,140)(340,212)
\qbezier(155,55)(300,200)(340,240)
\qbezier(170,235)(170,140)(170,56)
\qbezier(30,80)(100,140)(219,242) \qbezier(0,240)(100,140)(185,55)
\qbezier(0,240)(170,140)(320,60) \put(300,200){\circle*{3}}
\put(100,140){\circle*{3}}\put(170,200){\circle*{3}}\put(170,70){\circle*{3}}
\put(173,191){\small $P_3$}\put(303,191){\small $P_4$}
\put(96,127){\small $P_5$}\put(175,66){\small $P_6$}
\put(41,203){\small $P_1$}\put(159,149){\small $P_2$}
\put(220,203){\small $L_{(1,0,0)}$} \put(25,133){\small
$L_{(1,1,0)}$}\put(65,103){\small $L_{(0,0,1)}$}
\put(40,163){\small $L_{(1,1,1)}$}\put(138,230){\small
$L_{(1,0,1)}$}\put(290,80){\small $L_{(0,1,1)}$}
\put(250,140){\small $L_{(0,1,0)}$} \put(160,10){$\text{Fig. 6}$}
\end{picture}

Before, for convenience in further computations, we collect in
Table $(0,0,0,0)_2$ the adjacency types of the triangles
$P_1,\dots,P_6$ of $\mathcal{L}_{(0,0,0,0)}$. \vspace{0.2cm}

\begin{center}
{\small
\begin{tabular}{|c|c|c|c|c|} \hline
$P_1$ &   $(4_7,4_8',5_9,3_{15}',5_{14},4_{13}')$ & &  $P_2$ &
$(5_{16},3_{17}',5_{18},3_{21}',6_{20},3_{19}')$
\\ \hline
 $P_3$ &
$(5_9,4_{10}',5_{11},3_{17}',5_{16},3_{15}')$ & &
$P_4$ & $(5_{11},4_{12}',4_{13},4_7',5_{18}, 3_{17}')$ \\
\hline $P_5$ & $(4_{12},4_{13}',5_{14},3_{19}',6_{20},3_{22}')$ &
& $P_6$ & $(6_{20},3_{21}',4_{8},5_{9}',4_{10},3_{22}')$ \\
    \hline
\end{tabular}
} \end{center} \vspace{0.2cm}

\begin{center}
Table $(0,0,0,0)_2$
\end{center}
\vspace{0.3cm} (Here, we include in the adjacency type of the
triangle $P_i$, $i=1,\dots,6$, the indices of the adjacent
polygons. For example, in the adjacency type
$(4_7,4_8',5_9,3_{15}',5_{14},4_{13}')$ of $P_1$ the pattern $4_7$
points out that the polygon $P_7$ having four sides has a common
side with the triangle $P_1$.)

The adjacency type of the positive polygons of
$\mathcal{L}_{(0,0,0,0)}$ is equal to
$$A_{(0,0,0,0)}=((4,4',5,3',5,4'),(5,3',5,3',6,3')).$$

Perform in $\mathcal{L}_{(0,0,0,0)}$ a triangle move reversing
$P_3$. We obtain a new equipped purely real Campedelli line
arrangement. We denote it by $\mathcal{L}_{(1,0,0,0)}$ and we keep
to denote its polygons (denoted by $P'_i$ in subsection
\ref{reversing}) by $P_i$. To count its invariants, we notice,
first, that the adjacency type of $P_3$ changes as follows:
$$
(5_9,4_{10}',5_{11},3_{17}',5_{16},3_{15}')\mapsto
(4_9',5_{10},4_{11}',4_{17},4_{16}',4_{15}).
$$
After that, we adjust the number of sides of
$P_{9},P_{10},P_{11},P_{17},P_{16},P_{15}$ given in Tables
$(0,0,0,0)_1$ and $(0,0,0,0,)_2$ and obtain the tables for
$\mathcal{L}_{(1,0,0,0)}$: Table $(1,0,0,0)_1$ and Table
$(1,0,0,0)_2$.\vspace{0.3cm}

\begin{center}
\begin{tabular}{|c|c|c|c|c|c|c|c| 
} \hline
   $P_7$ & $P_8$ & $P_9$ & $P_{10}$ & $P_{11}$ & $P_{12}$ & $P_{13}$ & $P_{14}$ \\
   \hline
 4 & 4 & 4 & 5 & 4 & 4 & 4
& 5 \\ \hline \hline
 $P_{15}$
   & $P_{16}$ & $P_{17}$ & $P_{18}$ & $P_{19}$ & $P_{20}$ & $P_{21}$ &
   $P_{22}$
 \\ \hline 4 & 4 & 4 & 5 & 3 & 6 & 3 & 3 \\
   \hline
\end{tabular}
\end{center} \vspace{0.2cm}

\begin{center}
Table $(1,0,0,0)_1$
\end{center}
\vspace{0.2cm}

\begin{center}
{\small
\begin{tabular}{|c|c|c|c|c|} \hline
$P_1$ &   $(4_7,4_8',4_9,4_{15}',5_{14},4_{13}')$ & &  $P_2$ &
$(4_{16},4_{17}',5_{18},3_{21}',6_{20},3_{19}')$
\\ \hline
 $P_3$ &
$(4_9',5_{10},4_{11}',4_{17},4_{16}',4_{15})$ & &
$P_4$ & $(4_{11},4_{12}',4_{13},4_7',5_{18}, 4_{17}')$ \\
\hline $P_5$ & $(4_{12},4_{13}',5_{14},3_{19}',6_{20},3_{22}')$ &
& $P_6$ & $(6_{20},3_{21}',4_{8},4_{9}',5_{10},3_{22}')$ \\
    \hline
\end{tabular}
} \end{center} \vspace{0.2cm}

\begin{center}
Table $(1,0,0,0)_2$
\end{center}
\vspace{0.3cm}

We conclude that: $\mathcal{L}_{(1,0,0,0)}$ has the type
$(9,9,3,1,0)$; it contains two and only two positive polygons,
namely $P_1$ and $P_2$;
and the adjacency type of its positive
polygons is equal to
$$A_{(1,0,0,0)}=((4,4',4,4',5,4'),(4,4',5,3',6,3')).$$

Perform in $\mathcal{L}_{(1,0,0,0)}$ a triangle move reversing
$P_4$. Denote the new equipped purely real Campedelli line
arrangement by $\mathcal{L}_{(1,1,0,0)}$ and proceed as before. As
a result, we obtain two tables for $\mathcal{L}_{(1,1,0,0)}$:
Table $(1,1,0,0)_1$ and Table $(1,1,0,0)_2$.\vspace{0.3cm}

\begin{center}
\begin{tabular}{|c|c|c|c|c|c|c|c| 
} \hline
   $P_7$ & $P_8$ & $P_9$ & $P_{10}$ & $P_{11}$ & $P_{12}$ & $P_{13}$ & $P_{14}$ \\ \hline
 5 & 4 & 4 & 5 & 3 & 5 & 3
& 5 \\ \hline \hline $P_{15}$
   & $P_{16}$ & $P_{17}$ & $P_{18}$ & $P_{19}$ & $P_{20}$ & $P_{21}$ &
   $P_{22}$
 \\ \hline  4 & 4 & 5 & 4 & 3 & 6 & 3 & 3 \\
   \hline
\end{tabular}
\end{center}
\vspace{0.3cm}

\begin{center}
Table $(1,1,0,0)_1$
\end{center}
\vspace{0.4cm}

\begin{center}
{\small
\begin{tabular}{|c|c|c|c|c|} \hline

$P_1$ &   $(5_7,4_8',4_9,4_{15}',5_{14},3_{13}')$ & &  $P_2$ &
$(4_{16},5_{17}',4_{18},3_{21}',6_{20},3_{19}')$
\\ \hline
 $P_3$ &
$(4_9',5_{10},3_{11}',5_{17},4_{16}',4_{15})$ & &
$P_4$ & $(3_{11}',5_{12},3_{13}',5_7,4_{18}', 5_{17})$ \\
\hline $P_5$ & $(5_{12},3_{13}',5_{14},3_{19}',6_{20},3_{22}')$ &
& $P_6$ & $(6_{20},3_{21}',4_{8},4_{9}',5_{10},3_{22}')$ \\
    \hline
\end{tabular}
}\end{center} \vspace{0.3cm}

\begin{center}
Table $(1,1,0,0)_2$
\end{center}
\vspace{0.3cm}

We conclude that: $\mathcal{L}_{(1,1,0,0)}$ has the type
$(11,5,5,1,0)$; it contains two and only two positive polygons,
namely $P_1$ and $P_2$; and the adjacency type of its positive
polygons is equal to
$$A_{(1,1,0,0)}=((5,4',4,4',5,3'),(4,5',4,3',6,3')).$$

Perform in $\mathcal{L}_{(1,1,0,0)}$ a triangle move reversing
$P_3$. Denote the new equipped purely real Campedelli line
arrangement by $\mathcal{L}_{(0,1,0,0)}$ and proceed as before. As
a result, we obtain two tables for $\mathcal{L}_{(0,1,0,0)}$:
Table $(0,1,0,0)_1$ and Table $(0,1,0,0)_2$.\vspace{0.2cm}

\begin{center}
\begin{tabular}{|c|c|c|c|c|c|c|c| 
} \hline
   $P_7$ & $P_8$ & $P_9$ & $P_{10}$ & $P_{11}$ & $P_{12}$ & $P_{13}$ & $P_{14}$ \\ \hline
 5 & 4 & 5 & 4 & 4 & 5 & 3
& 5 \\ \hline \hline $P_{15}$
   & $P_{16}$ & $P_{17}$ & $P_{18}$ & $P_{19}$ & $P_{20}$ & $P_{21}$ &
   $P_{22}$
 \\ \hline  3 & 5 & 4 & 4 & 3 & 6 & 3 & 3 \\
   \hline
\end{tabular}
\end{center}
\vspace{0.2cm}

\begin{center}
Table $(0,1,0,0)_1$
\end{center}
\vspace{0.3cm}

\begin{center}
{\small
\begin{tabular}{|c|c|c|c|c|} \hline

$P_1$ &   $(5_7,4_8',5_9,3_{15}',5_{14},3_{13}')$ & &  $P_2$ &
$(5_{16},4_{17}',4_{18},3_{21}',6_{20},3_{19}')$
\\ \hline
 $P_3$ &
$(5_9,4_{10}',4_{11},4_{17}',5_{16},3_{15}')$ & &
$P_4$ & $(4_{11}',5_{12},3_{13}',5_7,4_{18}', 4_{17})$ \\
\hline $P_5$ & $(5_{12},3_{13}',5_{14},3_{19}',6_{20},3_{22}')$ &
& $P_6$ & $(6_{20},3_{21}',4_{8},5_{9}',4_{10},3_{22}')$ \\
    \hline
\end{tabular}
}\end{center} \vspace{0.2cm}

\begin{center}
Table $(0,1,0,0)_2$
\end{center}
\vspace{0.2cm}

We conclude that: $\mathcal{L}_{(0,1,0,0)}$ has the type
$(11,5,5,1,0)$; it contains two and only two positive polygons,
namely $P_1$ and $P_2$; and the adjacency type of its positive
polygons is equal to
$$A_{(0,1,0,0)}=((5,4',5,3',5,3'),(5,4',4,3',6,3')).$$

Perform in $\mathcal{L}_{(0,1,0,0)}$ a triangle move reversing
$P_5$. Denote the new equipped purely real Campedelli line
arrangement by $\mathcal{L}_{(0,1,1,0)}$ and proceed as before. As
a result, we obtain two tables for $\mathcal{L}_{(0,1,1,0)}$:
Table $(0,1,1,0)_1$ and Table $(0,1,1,0)_2$.\vspace{0.2cm}

\begin{center}
\begin{tabular}{|c|c|c|c|c|c|c|c| 
} \hline
   $P_7$ & $P_8$ & $P_9$ & $P_{10}$ & $P_{11}$ & $P_{12}$ & $P_{13}$ & $P_{14}$ \\ \hline
 5 & 4 & 5 & 4 & 4 & 4 & 4
& 4 \\ \hline \hline $P_{15}$
   & $P_{16}$ & $P_{17}$ & $P_{18}$ & $P_{19}$ & $P_{20}$ & $P_{21}$ &
   $P_{22}$
 \\ \hline  3 & 5 & 4 & 4 & 4 & 5 & 3 & 4 \\
   \hline
\end{tabular}
\end{center}
\vspace{0.2cm}

\begin{center}
Table $(0,1,1,0)_1$
\end{center}
\vspace{0.3cm}

\begin{center}
{\small
\begin{tabular}{|c|c|c|c|c|} \hline

$P_1$ &   $(5_7,4_8',5_9,3_{15}',4_{14},4_{13}')$ & &  $P_2$ &
$(5_{16},4_{17}',4_{18},3_{21}',5_{20},4_{19}')$
\\ \hline
 $P_3$ &
$(5_9,4_{10}',4_{11},4_{17}',5_{16},3_{15}')$ & & $P_4$ &
$(4_{11}',4_{12},4_{13}',5_7,4_{18}', 4_{17})$ \\
\hline $P_5$ & $(4_{12}',4_{13},4_{14}',4_{19},5_{20}',4_{22})$ &
 & $P_6$ & $(5_{20},3_{21}',4_{8},5_{9}',4_{10},4_{22}')$ \\
    \hline
\end{tabular}
}\end{center} \vspace{0.2cm}

\begin{center}
Table $(0,1,1,0)_2$
\end{center}
\vspace{0.2cm}

We conclude that: $\mathcal{L}_{(0,1,1,0)}$ has the type
$(8,10,4,0,0)$; it contains two and only two positive polygons,
namely $P_1$ and $P_2$; and the adjacency type of its positive
polygons is equal to
$$A_{(0,1,1,0)}=((5,4',5,3',4,4'),(5,4',4,3',5,4')).$$

Perform in $\mathcal{L}_{(0,1,1,0)}$ the triangle move reversing
$P_4$. Denote the new equipped purely real Campedelli line
arrangement by $\mathcal{L}_{(0,0,1,0)}$ and proceed as above. As
a result, we obtain two tables for $\mathcal{L}_{(0,0,1,0)}$:
Table $(0,0,1,0)_1$ and Table $(0,0,1,0)_2$.\vspace{0.2cm}

\begin{center}
\begin{tabular}{|c|c|c|c|c|c|c|c| 
} \hline
   $P_7$ & $P_8$ & $P_9$ & $P_{10}$ & $P_{11}$ & $P_{12}$ & $P_{13}$ & $P_{14}$ \\ \hline
 4 & 4 & 5 & 4 & 5 & 3 & 5
& 4 \\ \hline \hline $P_{15}$
   & $P_{16}$ & $P_{17}$ & $P_{18}$ & $P_{19}$ & $P_{20}$ & $P_{21}$ &
   $P_{22}$
 \\ \hline  3 & 5 & 3 & 5 & 4 & 5 & 3 & 4 \\
   \hline
\end{tabular}
\end{center}
\vspace{0.2cm}

\begin{center}
Table $(0,0,1,0)_1$
\end{center}
\vspace{0.3cm}

\begin{center}
{\small
\begin{tabular}{|c|c|c|c|c|} \hline

$P_1$ &   $(4_7,4_8',5_9,3_{15}',4_{14},5_{13}')$ & &  $P_2$ &
$(5_{16},3_{17}',5_{18},3_{21}',5_{20},4_{19}')$
\\ \hline
 $P_3$ &
$(5_9,4_{10}',5_{11},3_{17}',5_{16},3_{15}')$ & & $P_4$ &
$(5_{11},3_{12}',5_{13},4_7',5_{18}, 3_{17}')$ \\
\hline $P_5$ & $(3_{12}',5_{13},4_{14}',4_{19},5_{20}',4_{22})$ &
 & $P_6$ & $(5_{20},3_{21}',4_{8},5_{9}',4_{10},4_{22}')$ \\
    \hline
\end{tabular}
}\end{center} \vspace{0.2cm}

\begin{center}
Table $(0,0,1,0)_2$
\end{center}
\vspace{0.2cm}

We conclude that: $\mathcal{L}_{(0,0,1,0)}$ has the type
$(10,6,6,0,0)$; it contains two and only two positive polygons,
namely $P_1$ and $P_2$; and the adjacency type of its positive
polygons is equal to
$$A_{(0,0,1,0)}=((4,4',5,3',4,5'),(5,3',5,3',5,4')).$$

Perform in $\mathcal{L}_{(0,0,1,0)}$ a triangle move reversing
$P_3$. Denote the new equipped purely real Campedelli line
arrangement by $\mathcal{L}_{(1,0,1,0)}$ and proceed as before. As
a result, we obtain two tables for $\mathcal{L}_{(1,0,1,0)}$:
Table $(1,0,1,0)_1$ and Table $(1,0,1,0)_2$.\vspace{0.2cm}

\begin{center}
\begin{tabular}{|c|c|c|c|c|c|c|c| 
} \hline
   $P_7$ & $P_8$ & $P_9$ & $P_{10}$ & $P_{11}$ &
   $P_{12}$ & $P_{13}$ & $P_{14}$ \\ \hline
 4 & 4 & 4 & 5 & 4 & 3 & 5
& 4 \\ \hline \hline $P_{15}$
   & $P_{16}$ & $P_{17}$ & $P_{18}$ & $P_{19}$ & $P_{20}$ & $P_{21}$ &
   $P_{22}$
 \\ \hline  4 & 4 & 4 & 5 & 4 & 5 & 3 & 4 \\
   \hline
\end{tabular}
\end{center}
\vspace{0.2cm}

\begin{center}
Table $(1,0,1,0)_1$
\end{center}
\vspace{0.3cm}

\begin{center}
{\small
\begin{tabular}{|c|c|c|c|c|} \hline
$P_1$ &   $(4_7,4_8',4_9,4_{15}',4_{14},5_{13}')$ & &  $P_2$ &
$(4_{16},4_{17}',5_{18},3_{21}',5_{20},4_{19}')$
\\ \hline
 $P_3$ &
$(4_9',5_{10},4_{11}',4_{17},4_{16}',4_{15})$ & & $P_4$ &
$(4_{11},3_{12}',5_{13},4_7',5_{18}, 4_{17}')$ \\
\hline $P_5$ & $(3_{12}',5_{13},4_{14}',4_{19},5_{20}',4_{22})$ &
 & $P_6$ & $(5_{20},3_{21}',4_{8},4_{9}',5_{10},4_{22}')$ \\
    \hline
\end{tabular}
}\end{center} \vspace{0.2cm}

\begin{center}
Table $(1,0,1,0)_2$
\end{center}
\vspace{0.2cm}

We conclude that: $\mathcal{L}_{(1,0,1,0)}$ has the type
$(8,10,4,0,0)$; it contains two and only two positive polygons,
namely $P_1$ and $P_2$; and the adjacency type of its positive
polygons is equal to
$$A_{(1,0,1,0)}=((4,4',4,4',4,5'),(4,4',5,3',5,4')).$$

Finally, perform in $\mathcal{L}_{(1,0,1,0)}$ a triangle move
reversing $P_6$, denote the new equipped purely real Campedelli
line arrangement $\mathcal{L}_{(1,0,1,1)}$ and, once more, proceed
as before. As a result, we obtain two tables for
$\mathcal{L}_{(1,0,1,1)}$: Table $(1,0,1,1)_1$ and Table
$(1,0,1,1)_2$.\vspace{0.2cm}

\begin{center}
\begin{tabular}{|c|c|c|c|c|c|c|c| 
} \hline
   $P_7$ & $P_8$ & $P_9$ & $P_{10}$ & $P_{11}$ & $P_{12}$ & $P_{13}$ & $P_{14}$ \\ \hline
 4 & 3 & 5 & 4 & 4 & 3 & 5
& 4 \\ \hline \hline $P_{15}$
   & $P_{16}$ & $P_{17}$ & $P_{18}$ & $P_{19}$ & $P_{20}$ & $P_{21}$ &
   $P_{22}$
 \\ \hline  4 & 4 & 4 & 5 & 4 & 4 & 4 & 5 \\
   \hline
\end{tabular}
\end{center}
\vspace{0.2cm}

\begin{center}
Table $(1,0,1,1)_1$
\end{center}
\vspace{0.3cm}

\begin{center}
{\small
\begin{tabular}{|c|c|c|c|c|} \hline
$P_1$ &   $(4_7,3_8',5_9,4_{15}',4_{14},5_{13}')$ & &  $P_2$ &
$(4_{16},4_{17}',5_{18},4_{21}',4_{20},4_{19}')$
\\ \hline
 $P_3$ &
$(5_9',4_{10},4_{11}',4_{17},4_{16}',4_{15})$ & &
$P_4$ & $(4_{11},3_{12}',5_{13},4_7',5_{18}, 4_{17}')$ \\
\hline $P_5$ & $(3_{12}',5_{13},4_{14}',4_{19},4_{20}',5_{22})$ &
& $P_6$ & $(4_{20}',4_{21},3_{8}',5_{9},4_{10}',5_{22})$ \\
    \hline
\end{tabular}
}\end{center} \vspace{0.2cm}

\begin{center}
Table $(1,0,1,1)_2$
\end{center}
\vspace{0.2cm}

We conclude that: $\mathcal{L}_{(1,0,1,1)}$ has the type
$(8,10,4,0,0)$; it has two and only two positive polygons, namely
$P_1$ and $P_2$; and the adjacency type of its positive polygons
is equal to
$$A_{(1,0,1,1)}=((4,3',5,4',4,5'),(4,4',5,4',4,4')).$$

The results obtained show that each two of the eight constructed
arrangements either have different types or if their types
coincide, they have different adjacency types of their positive
polygons.

\subsection{Mixed real arrangements}\label{finmixreal}

\begin{ex}\label{ex-mixe}
Real Campedelli surfaces of types $I_{-}$ and
$III_{+}$ have diffeomorphic real structures, while they are not
deformation equivalent.
\end{ex}

Indeed, let $(X,c)$ be a real Campedelli surface of type $III_{+}$
associated with a mixed real Campedelli line arrangement of type
$III$. We can move the vertex $p_2$ so that it goes from the
triangle $P_3$ to $P_1$ through the line at infinity,
$L_{(1,1,0)}$. Theorem \ref{throughtriple} applies and shows that
the real structures $c$ and $c_1$ are diffeomorphic.

On the other hand, by Theorem \ref{def-c-c}, real surfaces of type
$I_-$ are not deformation equivalent to real surfaces of type
$III_+$.
\qed

In fact, in the case of mixed real types one can get a complete
answer to the Dif$\ne$Def. As it follows from the next theorem and
Theorem \ref{def-c-c}, the number of Dif classes is four, and the number of
Def classes is five.

\begin{thm} \label{dif-c-c}
The real structures of types
$I_{\pm}$, $II_+$, and
$III_{-}$ are pairwise non-diffeomorphic.
\end{thm}

\proof As it follows from Lemma \ref{ll}, their real points sets
have different topological types.\qed

\section{Final remarks}

\subsection{A pre-maximal surface}
 One can show that there are no $M$-
and ($M-1$)-surfaces among real Campedelli surfaces. As seems for
us, the following $(M-2)$-surface is of certain interest.

Let $\mathcal L$ be the purely real Campedelli line arrangement
depicted in Fig. 7. Its type is $(7,14,0,0,1)$ and it has three
and only three positive polygons: the $7$-gon $P_1$ and two
quadrangles, $P_2$ and $P_3$, with the sides $L_{(1,1,0)}$,
$L_{(0,1,1)}$, $L_{(0,0,1)}$, and $L_{(1,1,1)}$ for $P_2$, and
$L_{(0,1,0)}$, $L_{(1,0,0)}$, $L_{(0,1,1)}$, and $L_{(1,0,1)}$ for
$P_3$.

\begin{picture}(350,370)
\qbezier(20,231)(130,231)(370,231)
\qbezier(104,14)(231,172)(358,330)
\qbezier(35,190)(200,250)(365,310)
\qbezier(20,240)(150,170)(280,70)
\qbezier(50,96)(210,161)(370,236)
\qbezier(40,350)(125,182)(210,14)
\qbezier(20,350)(180,245)(340,140)\put(90,40){\small
$L_{(0,1,1)}$} \put(231,275){\small $L_{(1,0,1)}$}
\put(38,290){\small $L_{(1,1,0)}$} \put(220,235){\small
$L_{(1,0,0)}$} \put(270,80){\small $L_{(0,0,1)}$}
\put(50,112){\small $L_{(0,1,0)}$} \put(82,310){\small
$L_{(1,1,1)}$} \put(170,200){\small{$P_1$}}
\put(170,180){\small{$+++$}}\put(230,215){\small{$---$}}
\put(164,231){\small{$-++$}} \put(250,255){\small{$+--$}}
\put(320,255){\small{$+++$}} \put(340,270){\small{$P_3$}}
\put(90,271){\small{$++-$}}\put(170,300){\small{$--+$}}
\put(104,224){\small{$-+-$}} \put(63,218){\small{$+--$}}
\put(40,265){\small{$---$}} \put(23,210){\small{$+-+$}}
\put(84,202){\small{$--+$}}
\put(63,153){\small{$---$}} \put(143,144){\small{$++-$}}
\put(158,123){\small{$+--$}} \put(240,150){\small{$++-$}}
\put(179,142){\small{$+-+$}}\put(190,100){\small{$+++$}}
\put(200,75){\small{$P_{2}$}} \put(120,100){\small{$-+-$}}
\put(241,179){\small{$\times$}}
\qbezier(248,185)(249,186)(251,188)
\put(285,215){\small{$-++$}}\qbezier(253,190)(254,191)(256,193)
\put(295,185){\small{$--+$}} \put(160,-10){$\text{Fig. 7}$}
\end{picture}

\vspace{0.8cm}

Consider the Campedelli surface $X=X(\mathcal L)$ with its real
structure $c=c_{+++}$. As it follows from Corollary \ref{lll}, the
real part $X_{\mathbb R}$ of $(X,c)$ consists of three connected
components: the one over $P_1$ is a connected sum of eight real
projective planes (the Euler characteristic $-6$), the one over
$P_2$ is a Klein bottle, and the one over $P_3$ is a torus. It may
be interesting to notice that, in accordance with the Smith-Thom
inequality, $\dim H_*(X_\mathbb R;\mathbb Z/2\mathbb Z)=
10+4+4<22= \dim H_*(X_\mathbb C; \mathbb Z/2\mathbb Z)$, while the
ordinary Betti numbers of $X_{\mathbb R}$ surpass those of
$X_{\mathbb C}$: $\dim H_*(X_\mathbb R;\mathbb Q)= 8+2+4=14 > 10=
\dim H_*(X_\mathbb C; \mathbb Q)$.

\subsection{Bad moves}\label{badmove} Let us
show that the hypothesis on the signs in Corollary
\ref{rev-depend} is essential: without it, there is no local
equivariant diffeomorphism between the real Campedelli surfaces
$X(\mathcal L)$ and $X(\mathcal{L}')$ related by a triangle move
as in Corollary \ref{rev-depend}. For example, in the less evident
case, when $X(\mathcal L)$ (and thus $X(\mathcal L')$ as well) has
a real component over the triangle, to prove the nonexistence of a
local equivariant diffeomorphism one can argue in the following
way. We need to compare the quotients by the complex conjugation
of the Galois $(\mathbb Z/2\mathbb Z)^2$-coverings of a small ball
around the triple point ramified in, respectively, $ \mathcal L$
and $\mathcal L'$ (recall that $\alpha_1+\alpha_2+\alpha_3=0$).
More precisely, the boundaries of these quotients are naturally
identified, and the question is about possibility to extend this
identification to the interior. In fact, this is exactly one of
the questions treated in \cite{Pr} in an equivalent form, and as
it follows from \cite{Pr}, the extension {\it does not exist} if
and only if the four-manifold $M$ obtained by sewing of the
quotients along the boundary has the same homology as the
four-sphere $S^4$. Observing that there is a loop on the real
projective plane lying over the triangle in one half of $M$ linked
with the real projective plane lying over the triangle in the
other half, and using the Alexander-Pontryagin duality, one can
easily deduce that $H_*(M)=H_*(S^4)$.

\subsection{Class $T$ moves.} Theorem \ref{throughtriple}
(and its Corollary \ref{rev-depend}) are based on smoothings of
$T(-4)$-singularities. The latters constitute the simplest example of the
so-called {\it class $T$ singularities}. These singularities,
introduced by J.~Koll\'ar and N.~J.~Shepherd-Barron in
\cite{K-SB}, play crucial role in Manetti's examples \cite{Man}:
as is proved in \cite{Man}, smoothings of such singularities
provide diffeomorphic surfaces. As a consequence, the statement
and the proof of Theorem \ref{throughtriple} extend word-by-word
to real smoothings of real surfaces with any non real class $T$ singularities.

\subsection{On the number of deformation classes}\label{Finashinshina}
According to Proposition \ref{defeq}, the set of deformation
classes of real Campedelli surfaces splits into two disjoint
subsets: deformation classes of real surfaces associated with
mixed real Campedelli line arrangements, and those of real
surfaces associated with purely real Campedelli line arrangements.
By Theorem \ref{def-c-c}, the first subset contains only five
elements. Let us show that the other one contains more than
hundread elements.

We base our count on Proposition \ref{purerealclas} (and
Proposition \ref{defeq}), which imply that if $X(\mathcal L_1)$
and $X(\mathcal L_2)$ are real deformation equivalent, where
$\mathcal L_i$, $i=1,2$, are equipped purely real Campedelli line
arrangements without triple points, then, after a change of the
labels and the sign-equipment in $\mathcal L_2$ by a renumbering
homomorphism $h:({\mathbb Z}/2{\mathbb Z})^3\to ({\mathbb
Z}/2{\mathbb Z})^3$, it is possible to find a homeomorphism
$\lambda : \mathbb R\mathbb P^2\to \mathbb R\mathbb P^2$ which
transforms $\mathcal L_1\cap \mathbb R\mathbb P^2$ in $\mathcal
L_2\cap\mathbb R\mathbb P^2$ and preserves the labels and the
sign-equipments.

Consider an arrangement $\mathcal L$ of seven real lines which has
no triple points and is of type $(11,5,5,1,0)$. It has $7!$
distinct labelings turning it in a labelled purely real Campedelli
line arrangement, and, for each labeling, eight distinct
sign-equipments. Any homeomorphism of $\mathbb RP^2$ preserving
$\mathcal L\cap\mathbb R\mathbb P^2$ should preserve
$L_{(1,0,0)}\cap\mathbb R\mathbb P^2$ and the six-gon $P_{20}$
(see Fig. 5). It is easy to see that, up to isotopy fixing
$\mathcal L\cap\mathbb R\mathbb P^2$, there is only one such
homeomorphism, except identity. Since the order of the group
$\text{Aut} G$ of $G=(\mathbb Z/2\mathbb Z)^3$ is equal to $7\cdot
6\cdot 4=168$, we find that there are at least $\frac{7!\cdot
8}{(7\cdot 6\cdot 4)\cdot 2}=120$ distinct deformation classes of
real Campedelli surfaces $X=X(\mathcal L)$,  where $\mathcal L$ is
of the type $(11,5,5,1,0)$.

In fact, the number of deformation classes is even bigger. Indeed,
similar arguments show that at least $120$ more deformation
classes of real Campedelli surfaces are given by $X=X(\mathcal
L)$, where $\mathcal L$ is of the type $(9,9,3,1,0)$. In addition,
as is known (see \cite{Fi}, \cite{Cu}, and \cite{Wh}) there are
nine other deformation classes (of seven other types) of purely
real arrangements of seven lines without triple points. Note also
that two such arrangements are deformation equivalent if, and only
if, they are homeomorphic, see \cite{Fi} (proofs are found
in \cite{FiT}). Similarly, a self-homeomorphism of
an equipped purely real arrangement of seven lines without triple
points should be isotopic, together with the arrangement, to a
projective automorphism, which would imply, according to Corollary
\ref{defeq}, that the number of deformation classes of purely real
Campedelli surfaces is the same as the number of deformation
classes of equipped purely real Campedelli arrangements without
triple points.

\ifx\undefined\bysame
\newcommand{\bysame}{\leavevmode\hbox to3em{\hrulefill}\,}
\fi

\end{document}